\documentclass[a4]{article}%For paper submission

\pdfoutput=1

\usepackage{graphicx}
\usepackage{times}
\usepackage{pifont,latexsym,ifthen,rotating,calc,textcase,booktabs,color}
\usepackage{amsfonts,amssymb,amsbsy,amsmath,amsthm}

\usepackage{xcolor}
\usepackage[ruled,section]{algorithm}
\usepackage{algorithmic}
\usepackage{multirow}

\newcommand{\gradot}{{}^{t}\underline{\underline{\nabla}}_\textit0}
\newcommand{\grado}{\underline{\underline{\nabla}}_\textit0}
\newcommand{\x}{{\underline{x}}}
\newcommand{\X}{{\underline{X}}}
\newcommand{\ply}{P}
\newcommand{\plyp}{P'}
\newcommand{\plyo}{{P^{\phantom{'}}_\textit0}}
\newcommand{\plypo}{{P'_\textit0}}

\newcommand{\uply}{{\underline{u}_{PP'}}}

\newcommand{\tracply}{{\underline{t}_P}}
\newcommand{\tracplyo}{{\underline{t}_{P_\textit0}}}

\newcommand{\FFply}{\underline{\underline{F}}_{PP'}}

\newcommand{\GammaPPo}{{\Gamma_{P^{\phantom{'}}_\textit0P'_\textit0}}}
\newcommand{\GammaPP}{{\Gamma_{PP'}}}

\newcommand{\Fd}{{\underline{F}_d}}
\newcommand{\Ud}{{\underline{U}_d}}
\newcommand{\fd}{{\underline{f}_d}}
\newcommand{\normalply}{{\underline{n}_3}}
\newcommand{\tanplyone}{{\underline{n}_1}}
\newcommand{\tanplytwo}{{\underline{n}_2}}
\newcommand{\normalplyo}{{\underline{N}_{3}}}
\newcommand{\Kinter}{{\underline{\underline{K}}_{PP'}}}
\newcommand{\rota}{{\underline{\underline{Q}}}}
\newcommand{\umoyPP}{{{\langle\underline{u}_{PP'}\rangle}}}
\newcommand{\Confio}{\mathcal{C}_\textit0}
\newcommand{\struct}{{\mathbf{E}}}

\newcommand{\Omegao}{{\Omega_\textit0}}

\newcommand{\dOmegaUd}{{\partial \Omega_{U_d}}}
\newcommand{\dOmegaFd}{{\partial \Omega_{F_d}}}
\newcommand{\dOmegaUdo}{{\partial \Omega_{U_{d_\textit0}}}}
\newcommand{\dOmegaFdo}{{\partial \Omega_{F_{d_\textit0}}}}
\newcommand{\Fdo}{{\underline{F}_{d_\textit0}}}
\newcommand{\Udo}{{\underline{U}_{d_\textit0}}}
\newcommand{\GammaEdo}{{\Gamma_{{E}_{d_\textit0}}}}
\newcommand{\GammaEo}{{\Gamma_{E_\textit0}}}
\newcommand{\GammaEEpo}{{\Gamma_{E^{\phantom{'}}_\textit0E'_\textit0}}}

\newcommand{\Gammao}{{\Gamma_\textit0}}
\newcommand{\E}{{E^{\phantom{'}}_\textit0}}
\newcommand{\Ep}{{E'_\textit0}}
\newcommand{\OmegaEo}{{\Omega_{E_\textit0}}}
\newcommand{\OmegaEpo}{{\Omega_{E_\textit0'}}}
\newcommand{\dOmegaEo}{{\partial \Omega_{E_\textit0}}}
\newcommand{\uE}{{\underline{u}_{E^{\phantom'}_\textit0}}}

\newcommand{\uEstar}{{\underline{u}^{\star}_{E_\textit0}}}

\newcommand{\EpsilonpE}{{\underline{\underline{\dot{E}}}}}
\newcommand{\EpsilonpEstar}{{\underline{\underline{\dot{E}}}(\underline{u}^{\star}_{E_\textit0})}}
\newcommand{\GLE}{{\underline{\underline{E}}_{E_\textit0}}}
\newcommand{\GLpE}{{\underline{\underline{\dot{E}}}_{E_\textit0}}}
\newcommand{\piE}{{\underline{\underline{\pi}}_{E_\textit0}}}
\newcommand{\KEo}{\mathbf{K}_{E_\textit0}}

\newcommand{\FF}{\underline{\underline{F}}}
\newcommand{\FFE}{\underline{\underline{F}}_{E^{\phantom{'}}_\textit0}}
\newcommand{\FFEE}{\underline{\underline{F}}_{EE'}}

\newcommand{\WmoyEE}{{{\langle\underline{W}_{E^{\phantom'}_\textit0}\rangle}}}
\newcommand{\WE}{{\underline{W}_{E^{\phantom'}_\textit0}}}
\newcommand{\WEM}{{\underline{W}_{E_\textit0}^M}}
\newcommand{\WEm}{{\underline{W}_{E_\textit0}^m}}
\newcommand{\WEstar}{{{\underline{W}^{\star}_{E_\textit0}}}}
\newcommand{\FEstar}{{{\underline{F}^{\star}_{E_\textit0}}}}
\newcommand{\WMstar}{{{\underline{W}^{M}}}^\star}
\newcommand{\Wmstar}{{\underline{W}^{m}}^\star}
\newcommand{\WEp}{{\underline{W}_{E'_\textit0}}}
\newcommand{\FEo}{{\underline{F}_{E^{\phantom'}_\textit0}}}

\newcommand{\FEoM}{{\underline{F}^M_{E^{\phantom'}_\textit0}}}
\newcommand{\FEom}{{\underline{F}^m_{E^{\phantom'}_\textit0}}}
\newcommand{\FEpo}{{\underline{F}_{E_\textit0'}}}

\newcommand{\kmEo}{{k^-_{E^{\phantom'}_\textit0}}}

\newcommand{\kmEoM}{{k^{-M}_{E^{\phantom'}_\textit0}}}
\newcommand{\kmEom}{{k^{-m}_{E^{\phantom'}_\textit0}}}

\newcommand{\kpEo}{{k^+_{E^{\phantom'}_\textit0}}}

\newcommand{\kpEpo}{{k^+_{E_\textit0'}}}
\newcommand{\sWE}{{\mathcal{W}_E}}
\newcommand{\sWEO}{{\mathcal{W}_E^0}}

\newcommand{\sFE}{{\mathcal{F}_{E_\textit0}}}
\newcommand{\sFEO}{{\mathcal{F}_{E}^\textit0}}

\newcommand{\sAd}{\mathbf{A_d}}
\newcommand{\sGamma}{{\boldsymbol{\Gamma}}}

\newcommand{\suEO}{{\mathcal{U}_E^0}}

\newcommand{\sFM}{{\mathcal{F}^M}}
\newcommand{\sWM}{{\mathcal{W}^M}}

\newcommand{\sWm}{{\mathcal{W}^m}}

\newcommand{\sEp}{{\mathbf{E^+}}}
\newcommand{\sEm}{{\mathbf{E^-}}}
\newcommand{\sWMad}{{\mathcal{W}^M_\textrm{ad}}}

\newcommand{\FFchap}{{\underline{\underline{\widehat{F}}}}}
\newcommand{\FchapEo}{{\underline{\widehat{F}}_{E^{\phantom'}_\textit0}}}
\newcommand{\FchapEpo}{{\underline{\widehat{F}}_{E'_\textit0}}}
\newcommand{\FchapchapEo}{{\underline{\widehat{\widehat{F}}}_{E_\textit0}}}

\newcommand{\WchapE}{{\underline{\widehat{W}}_{E^{\phantom'}_\textit0}}}
\newcommand{\WchapEp}{{\underline{\widehat{W}}_{E'_\textit0}}}

\newcommand{\Wtilde}{{{}^{i+1}\underline{\widetilde{W}}^M}}
\newcommand{\Wtildestar}{{{\underline{\widetilde{W}}^M}^{\star}}}
\newcommand{\WtildeE}{{}^{i+1}\underline{\widetilde{W}}_{E^{\phantom'}_\textit0}^M}
\newcommand{\FtildeEo}{{}^{i+1}{\underline{\widetilde{F}}_{E^{\phantom'}_\textit0}^M}}

\newcommand{\WEdc}{{}^i\underline{W}_{E^{\phantom'}_\textit0}^c}
\newcommand{\LEM}{{}^{i}\mathbb{L}_{E^{\phantom'}_\textit0}^M}
\newcommand{\HE}{{}^i\mathbb{H}_{E^{\phantom'}_\textit0}}

\graphicspath{{figures/}}

\begin{document}
\markboth{K.~Saavedra, O.~Allix and P.~Gosselet}{On a multiscale strategy and its optimization for delamination and buckling}
\title{On a multiscale strategy and its optimization for the simulation of combined delamination and buckling}
\author{K. Saavedra, O. Allix and P. Gosselet\\
LMT Cachan (ENS Cachan/CNRS/UPMC/PRES UniverSud Paris) \\61,~Avenue~du~Pr{\'e}sident~Wilson,~94235~Cachan,~France.\\
E-mail: [saavedra,allix,gosselet]@lmt.ens-cachan.fr}
%\cgs{<Contract/grant sponsor name (no number)>}
%\cgsn{<Contract/grant sponsor name>}{<number>}
\date{2012}

\maketitle

\begin{abstract}
This paper investigates a computational strategy for studying the interactions between multiple through-the-width delaminations and global or local buckling in composite laminates taking into account possible contact between the delaminated surfaces. In order to achieve an accurate prediction of the quasi-static response, a very refined discretization of the structure is required, leading to the resolution of very large and highly nonlinear numerical problems. In this paper, a nonlinear finite element formulation along with a parallel iterative scheme based on a multiscale domain decomposition are used for the computation of 3D mesoscale models. {Previous works by the authors already dealt with the simulation of multiscale delamination assuming small perturbations. This paper presents the formulation used to include geometric nonlinearities into this existing} multiscale framework and discusses the adaptations that need to be made to the iterative process in order to ensure the rapid convergence and the scalability of the method in the presence of buckling and delamination. These various adaptations are illustrated by simulations involving large numbers of DOFs.
\end{abstract}
%\medskip

\textbf{keywords}: nonlinear multiscale computation; domain decomposition method; delamination; buckling; composites
%\medskip 

\section{Introduction}
\label{sec:intro}

Delamination is one of the main degradation mechanisms of laminated composite materials. This phenomenon is generally initiated by large interlaminar stresses due to edge effects, impacts, concentrated loads or macroscopic defects. Under some loading and geometric configurations (\emph{e.g.}  compression and a high slenderness coefficient), buckling is likely to occur once the delaminated zone has reached a critical extent during the propagation phase of the delamination process. Then, this geometric instability can lead to an increase in interlaminar stresses, an acceleration of the delamination rate and, eventually, to the failure of the structure. The first analytical studies of buckling and delamination growth in the 70-80's were done by \cite{Kachanov76, Chai81, Bottega83, Evans84}; in the last decade, new analytical studies were proposed by \cite{Daridon02, Pradeilles-Duval04, Kosel05}. In a finite element context, the first works were based on fracture mechanics \cite{Storakers88, Whitcomb89, Nilsson90} while in more recent publications cohesive models were also used to deal with geometrically nonlinear problems \cite{Bruno90,Allix99,Qiu01}. Asymptotic numerical methods \cite{Cochelin91,Kardomateas93, Bruno00} were also applied to delamination buckling problems. Despite these many contributions to a better understanding of the mechanics of laminated composites, the inclination by industry to substitute virtual simulations for expensive experimental tests raises new issues. Thus, the numerical prediction of combined buckling and delamination remains a scientific challenge because, even when using calculations on the mesoscale \cite{Ladeveze02}, a highly refined discretization of each ply is necessary in order to describe the delamination fronts and buckling loads properly. Therefore, multiscale and parallel computational techniques are being developed for buckling \cite{Cresta07,Nezamabadi10} and debonding problems \cite{Guidault08,Kerfriden09}.

{In this work, we propose a mixed and multiscale domain decomposition strategy for the parallel simulation, in a geometrically nonlinear context, of composite laminates which are subject to multiple delaminations. Our approach was adapted to the treatment of geometric nonlinearities from an existing LATIN (LArge Time INcrement) multiscale strategy for delamination problems under the assumption of small perturbations \cite{Kerfriden09}.}

Here, the geometrically nonlinear evolution is handled through a total Lagrangian formulation and delamination is modeled on the mesoscale using a cohesive interface model based on damage mechanics \cite{Allix99}. For this first-time approach, the intralaminar degradations are considered to be negligible and the layers are assumed to follow a hyperelastic law. Unilateral contact conditions are introduced by means of an interface law in order to avoid interpenetration over the delaminated surfaces. The reference problem and the substructuring process are summarized in Section \ref{section:ref_pb}.

The LATIN strategy {\cite{Ladeveze99}} consists in dividing the structure into volume substructures separated by 2D interfaces, both of which are mechanical entities. As a result, the reference problem associated with the chosen mesomodel is naturally substructured, and both the unilateral contact conditions and the cohesive interfaces are handled at the interfaces of the domain decomposition. {Section \ref{section:strategy} introduces the LATIN algorithm proposed for the resolution of the nonlinear substructured problem.} For the sake of computational efficiency, three scales are considered in this resolution:
\begin{itemize}
  \item The \textbf{microscale} corresponds to small-wavelength phenomena, which occur between neighboring substructures.
  \item The \textbf{macroscale} corresponds to the permanent verification of a weak form of equilibrium throughout the structure. This part of the algorithm, which makes the method scalable, is achieved through the definition of a small number of macroscopic degrees of freedom per interface, which must satisfy continuity conditions and are linked together by a homogenized behavior constructed automatically.
  \item In some cases, the number of substructures and interfaces (which depends on the number of plies) may be so large that the macroscopic problem cannot be addressed by direct solvers. Therefore, the substructures are grouped into ``supersubstructures'' (whose size is determined by the available processor memory), and the macroscopic problem is solved using a primal domain decomposition method \cite{Mandel93}. The \textbf{third scale} (or supermacroscopic problem) is introduced classically in the course of balancing the supersubstructures with respect to the rigid body modes.
\end{itemize}

This framework is very much under the control of the operator. In order to pilot the calculation of slender structures, the following parameters must be adjusted:
\begin{itemize}
  \item The influence of neighboring subdomains and interfaces is represented by what is known as ``search directions'' (also called ``interface impedances'' in \cite{Ladeveze01}). These parameters must be adapted to the aspect ratios of the slender structures by introducing well-chosen anisotropic coefficients; this point is discussed in Section \ref{section:ddr_perfect_interfaces}.
  \item Unilateral contact with or without friction between small surfaces can be handled successfully by the multiscale LATIN method \cite{Champaney99, Ladeveze02a}. However, as shown in Section \ref{section:ddr_contact}, in the case of contact between slender structures over large delaminated areas, the search directions should be optimized according to the interface's state (open or closed) because incorrect values could generate artificial stiffnesses or induce interpenetration of the contact surfaces.
  \item Because of the stiffness loss as a result of buckling and delamination, the macrostiffness and search directions may become irrelevant and the macroscopic operators may need to be adjusted in order for the homogenized behavior to represent the current state of the structure better, as illustrated in Sections \ref{section:buckling} and \ref{section:ddr_cohesive_interfaces}.
  \item The supermacroscopic problem is solved using a projected preconditioned conjugate gradient algorithm, which requires the setting up of a convergence threshold. The solutions to this problem were developed in \cite{Kerfriden09} and will not be repeated here.
\end{itemize}

%%%%%%%%%%%
With these improvements, the multiscale analysis of large combined buckling and delamination problems becomes possible. The capabilities of the strategy are illustrated by two examples involving geometric instabilities coupled with debonding (Section \ref{section:flam_delam}).

\section{The reference problem}\label{section:ref_pb}

\subsection{Notations and assumptions}
\label{subsec:desc_motion}
\begin{figure}[ht]
       \centering
       \includegraphics[width=0.65 \linewidth]{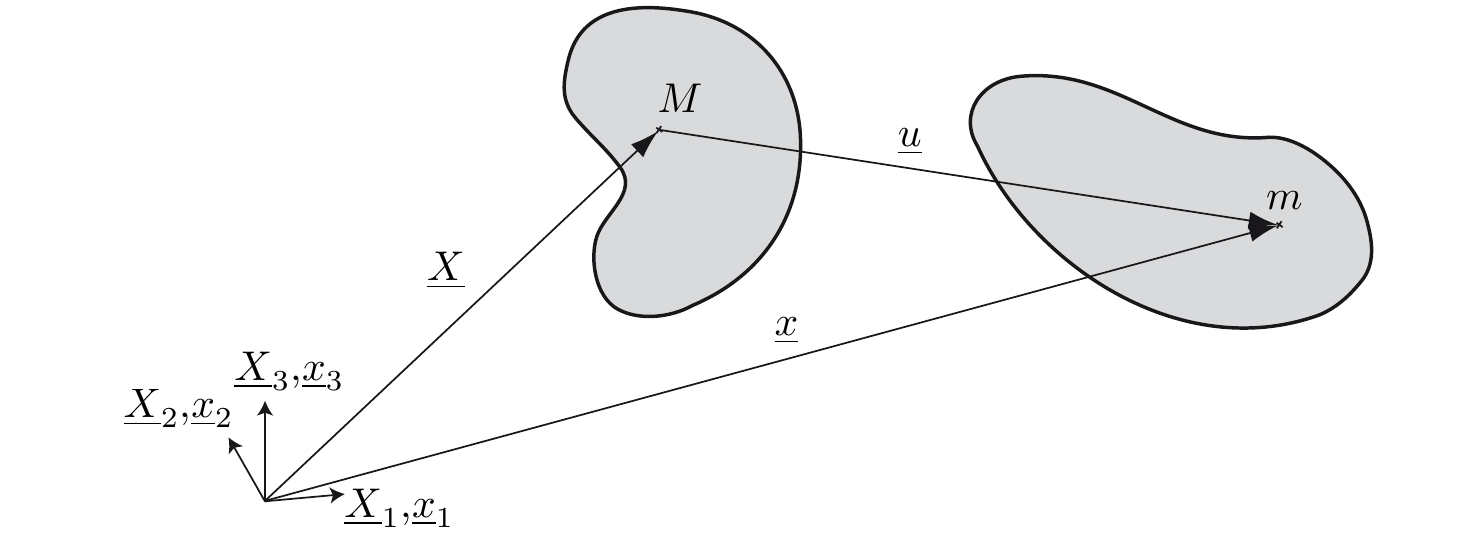}
       \caption{The general motion of a deformable body}
       \label{fig:motion}
\end{figure}

Figure \ref{fig:motion} shows the general motion of a deformable body. The body is considered to be an assembly of material particles $M$ identified by their initial coordinates $\X$ with respect to the Cartesian frame $\mathcal{B}_\textit0 = \{\X_1, \X_2, \X_3\}$. In general, the current positions of these particles are defined by their coordinates $\x$ with respect to another Cartesian frame $\mathcal{B} = \{\x_1, \x_2, \x_3\}$. In this work, the two coordinate systems $\mathcal{B}_\textit0$ and $\mathcal{B}$ are the same, but we will refer to them as separate entities in order to associate each quantity with the initial or the current configuration.
The displacement $\underline{u}$ of particle $M$ between the two configurations is defined as:
\begin{equation}
\underline{u} = \x - \X\;.\;
\end{equation}
The deformation gradient tensor $\FF$ is given by:
\begin{equation}
\FF = \grado \x = \underline{\underline{I}}_{d} + \grado \underline{u}\;,\;
\end{equation}
where $\grado$ denotes the gradient with respect to the initial configuration. Let us note that $\FF$, which transforms vectors in the initial configuration into vectors in the current configuration, is called a two-point tensor.

\begin{figure}[b]
       \centering
       \includegraphics[width=.65 \linewidth]{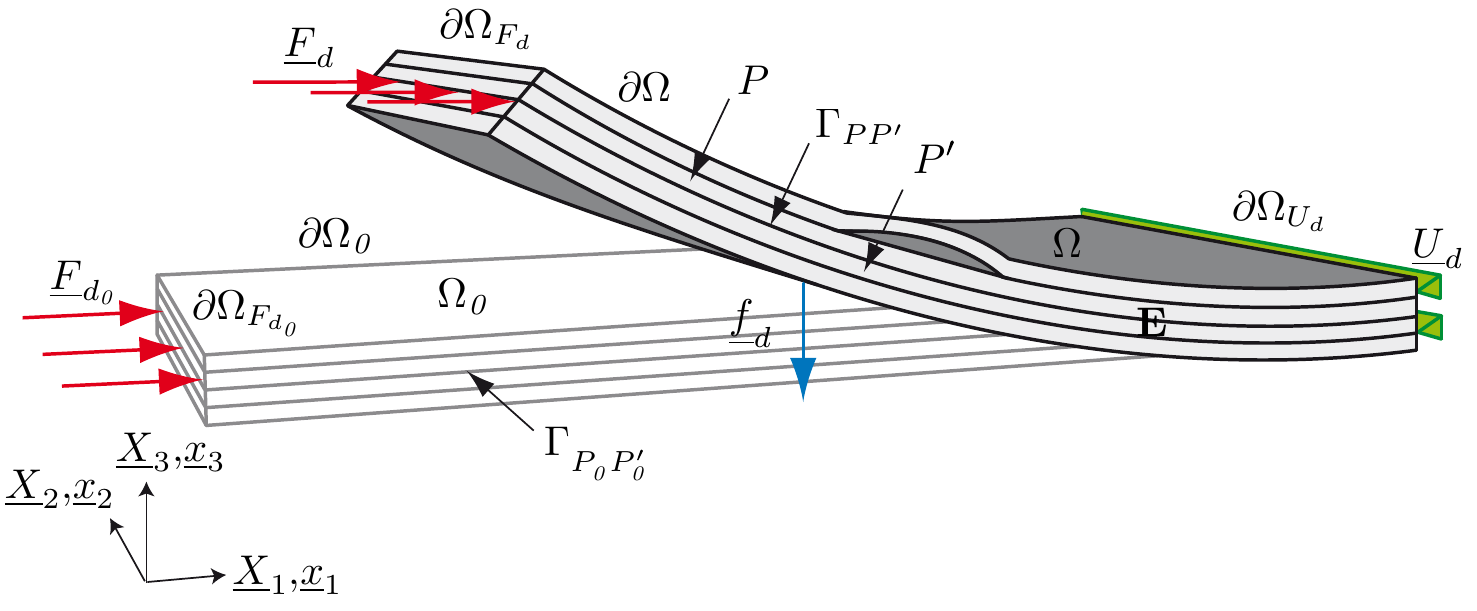}
       \caption{The reference problem: undeformed and deformed configurations}
       \label{fig:reference_pb}
\end{figure}

Let us consider a laminated composite structure $\struct$ (see Figure \ref{fig:reference_pb}) occupying domain $\Omega$ bounded by $\partial \Omega$ in the current configuration $\mathcal C$, and consisting of $N_P$ adjacent plies $\ply$. Each ply $\ply$, of mass density $\rho_P$, occupies a domain $\Omega_P$ such that $\Omega = \bigcup_{P \in \struct} \Omega_P$. The plies are assumed to be separated by $N_P-1$ cohesive interfaces. The structure is subjected to an external surface traction field $\Fd$ over part $\dOmegaFd$ of the boundary $\partial \Omega$ and to a displacement field $\Ud$) over the complementary part $\dOmegaUd$. The body force per unit mass is denoted $\fd$.
A ply $\ply$ defined in domain $\Omega_P$ is connected to an adjacent ply $\ply'$ through an interface $\GammaPP=\partial \Omega_P \cap \partial \Omega_{P'}$. Let $\Gamma = \bigcup_{P \in \struct} \Gamma_P$, where $\Gamma_P = \bigcup_{P' \in \struct} \GammaPP$. The relevant quantities of $\struct$ (\emph{e.g.}  volume, area, surface tractions or density) can be described in reference to the configuration before deformation $\Confio$. The index $\cdot_\textit0$ will be used to denote the initial (undeformed) configuration, \emph{e.g.}  $\Omegao$, $\partial \Omegao$, $ \GammaPPo$, $\Fdo$, $\rho_{P_\textit0}$.

The objective of the present work is to study the response of $\struct$ subjected to a prescribed loading starting from the initial configuration, and resulting in large displacements and rotations accompanied by progressive damage of the interfaces $\Gamma$. Because of the very small thickness of each layer of the composite ($\approx 0.125$ mm), the delaminated areas can be very slender $(L_{delaminated}/h_{ply} \gg 100)$.

Our study of this problem relies on the following assumptions:
\begin{enumerate}
\item structure $\struct$ may undergo large displacements;
\item the behavior of the plies is hyperelastic;
\item the loads are independent of the configuration of $\struct$ (\emph{i.e.} follower forces are not considered);
\item the evolution over time is considered to be quasi-static;
\item isothermal conditions are assumed;
\item the displacements along $\Gamma$ can be large, but the displacement discontinuities in the non-fully delaminated part of the structure or in the contact region are small; \label{hyp_inter}
\item $\Gamma$ is assigned irreversible softening behavior by means of an interface law connecting tractions with displacement discontinuities;
\item in the delaminated region, corresponding points in adjacent plies may separate, regain contact after separation, or remain in contact.
\end{enumerate}

Assumption \ref{hyp_inter} enables the displacement of interface $\GammaPP$ to be defined as the mean value of the displacements of the adjacent plies (see Figure \ref{fig:normal_interface}):
\begin{equation}
\umoyPP = \frac{1}{2}(\underline{u}_{P'}+\underline u_P) \;,\; \text{over} \; \GammaPP \;.
\end{equation}
Thus, interface $\GammaPP$ itself is defined as:
\begin{equation}
\GammaPP \equiv \{ \underline{x} \backslash \underline{x} = \underline{X} +  \umoyPP \,;\, \text{over} \; \GammaPP\} \;.
\end{equation}
With this definition, $\GammaPP$ is viewed as the mean surface between deformed plies $\ply$ and $\plyp$, which, in fact, due to Assumption \ref{hyp_inter}, can be considered to coincide geometrically (\emph{i.e.} same area and same orientation). Because of this geometric coincidence, the interface's deformation gradient tensor can be defined as:
\begin{equation}
\FF_{PP'} = \grado \umoyPP+ \underline{\underline{I}}_{d}  \;,\; \text{over} \; \GammaPP \;.
\end{equation}

Thus, interface $\GammaPP$ can be viewed as a zero-thickness medium that carries out the transfer of traction forces between plies. The displacement gap of interface $\GammaPP$ is given by:
\begin{equation}
[\uply] =\underline{u}_{P'}-\underline{u}_{P} \;,\; \text{over} \; \GammaPP\;.
\end{equation}

\subsection{The nonlinear damage interface law}
\label{section:interface_law}
An anisotropic interface law can be formulated as a relation between the traction vector $\tracply$ and the displacement discontinuity vector $[\uply]$ :  
\begin{equation}
\tracply = \Kinter([\uply]) \, [\uply] \;,\; \text{over} \; \GammaPP\;,
\label{eq:inter_law2}
\end{equation}
in which a dependence on the current orientation of the interface $\GammaPP$ is introduced. One can show, as proven in \cite{Allix99}, that law \eqref{eq:inter_law2} satisfies the principle of material frame indifference if Assumption \ref{hyp_inter} is satisfied.

\begin{figure}[ht]
       \centering
       \includegraphics[width=.65 \linewidth]{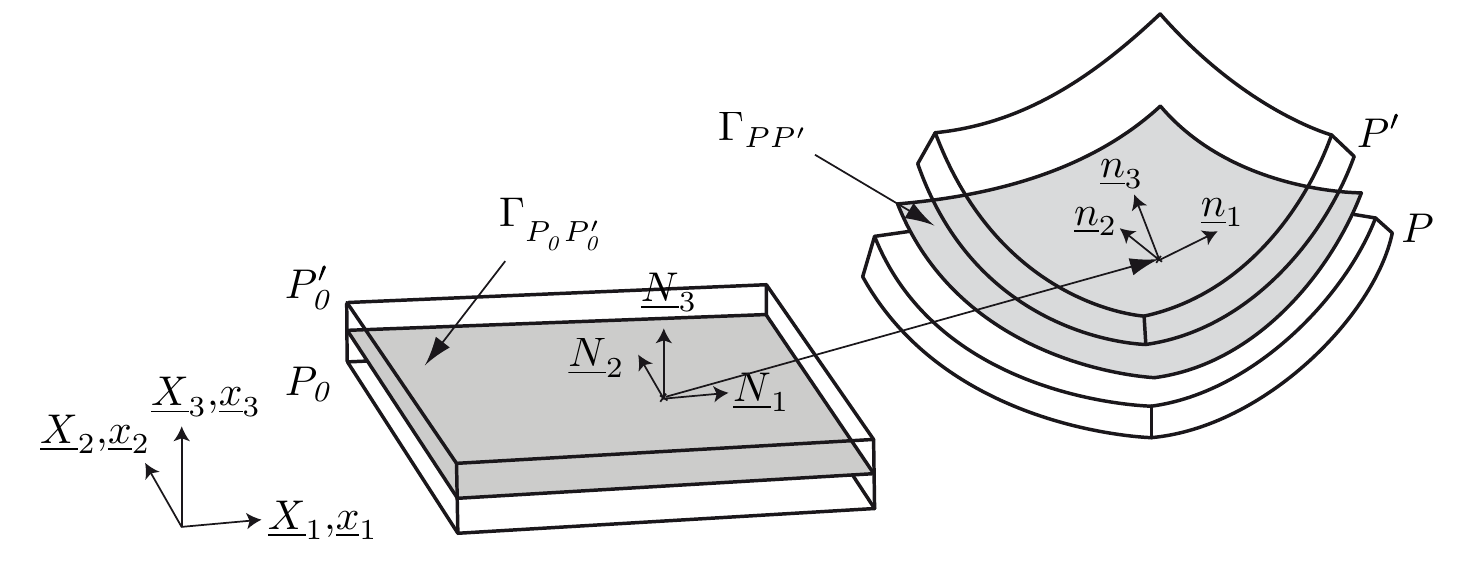}
       \caption{The normals to the undeformed and deformed interfaces}
       \label{fig:normal_interface}
\end{figure}

The expression of the local stiffness operator $\Kinter$ of interface $\GammaPP$ can be made explicit with respect to a local orthonormal frame $\mathcal{B}_n = \{\tanplyone,\tanplytwo,\normalply\}$ moving together with the interface, where $\normalply$ is the deformed unit normal to $\GammaPP$ pointing from $\ply$ to $\plyp$ (see Figure \ref{fig:normal_interface}):
\begin{equation}
\Kinter = K^{local}_{ij} \underline{n}_i \otimes \underline{n}_j \;,\;
\end{equation}
where:
\begin{equation*}
\nonumber
(K^{local}_{ij}) = \left( \begin{array}{ccc}
 \displaystyle (1-d_1)\, k_t^0 & 0 & 0 \\
0 & \displaystyle (1-d_2)\, k_t^0 & 0 \\
0 & 0 & \displaystyle \left(1-h_+([\uply] \cdot \normalply)\, d_3 \right) k_n^0
\end{array}
\right) \;,\;
\end{equation*}
$h_+$ denotes the positive indicator function. $k_n^0$ and $k_t^0$ are the initial elastic stiffnesses of the interface,
with the dimension of a force per volume. The softening behavior of the interface model when the structure is loaded is simulated by the introduction of the dimensionless scalar damage variables ${d_i}$ with values ranging from $0$ (healthy interface point) to $1$ (completely damaged interface point). {We use the evolution law defined in \cite{Allix98}, which has the advantage of using a single damage variable to handle different macroscopic delamination modes of the interface.}

To calculate the internal power corresponding to the interfaces, it is mandatory to express the traction vector as a function of the displacement discontinuity vector with respect to frame $\mathcal{B}_\textit0$. In order to do that, one needs to write $\Kinter = K^{global}_{ij} \underline{X}_i \otimes \underline{X}_j$. Introducing the orthogonal matrix which characterizes the transition from $\mathcal{B}_n$ to $\mathcal{B}_\textit0$:
\begin{equation}
\rota = Q_{ij} \underline{X}_i \otimes \underline{n}_j \quad \forall i,j = \{1,2,3\} \;,\;
\end{equation}
where $Q_{ij}= \underline{X}_i \cdot \underline{n}_j$ are the direction cosines of $\underline{X}_i$ relative to $\underline{n}_j$, the interface's local stiffness becomes:
\begin{equation}
(K^{global}_{ij}) = Q_{ip} K^{local}_{pq} Q_{jq} \quad \forall i,j = \{1,2,3\} \; \forall p,q = \{1,2,3\}
\end{equation}
(with implied summation over indexes $p,q$).

From Nanson's formula \footnote{$\normalply d\GammaPP = det(\FFply)\FFply^{-t} \normalplyo d\GammaPPo $}, which describes how an infinitesimal surface element deforms during a given motion, the unit normal vector at each point of the deformed interface $\GammaPP$ is calculated as follows:
\begin{equation}
\normalply = \frac{\FFply^{-t} \, \normalplyo}{\Vert \FFply^{-t} \, \normalplyo \Vert } \;,\;
\label{eq:deformed_normal}
\end{equation}
where $\normalplyo$ is the initial (undeformed) unit normal to $\GammaPPo$ pointing from $\plyo$ to $\plypo$ (see Figure \ref{fig:normal_interface}).

Finally, one needs to calculate the traction vector in the initial configuration:
\begin{equation}
\tracplyo = \frac{d\GammaPP}{d\GammaPPo}\, \tracply = det(\FFply) \, \Vert \FFply^{-t} \, \normalplyo \Vert \, \tracply \;.\;
\end{equation}

After an interface has been fully damaged, a frictionless contact law is assumed: the gap between two plies must remain nonnegative, and only compression can be transmitted when the plies are in contact.

\textit{Remark:} The relative interface displacements are assumed to remain small enough (Assumption \ref{hyp_inter}) for the contact to be detected only between points which were connected in the initial configuration.

\subsection{Substructured formulation}
\label{subsec:substructure}

\begin{figure}[tbh]
       \centering
       \includegraphics[width=.65 \linewidth]{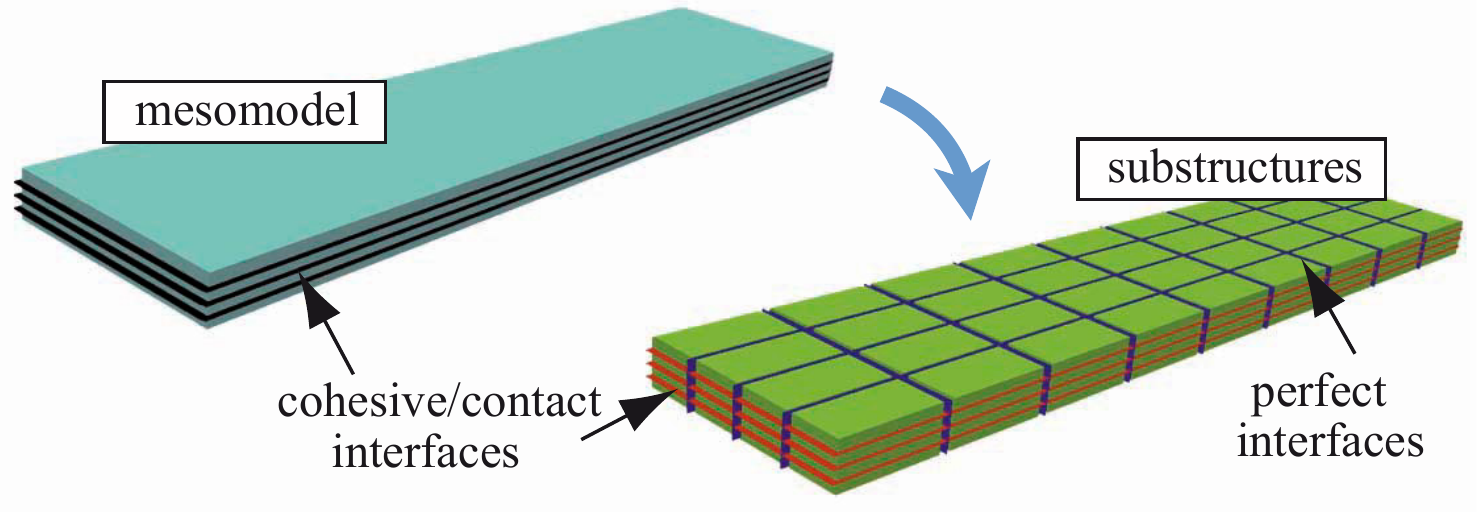}
       \caption{Substructuring of the laminated composite structure}
       \label{fig:decomp_sst_interfaces}
\end{figure}

\begin{figure}[htb]
       \centering
       \includegraphics[width=.65 \linewidth]{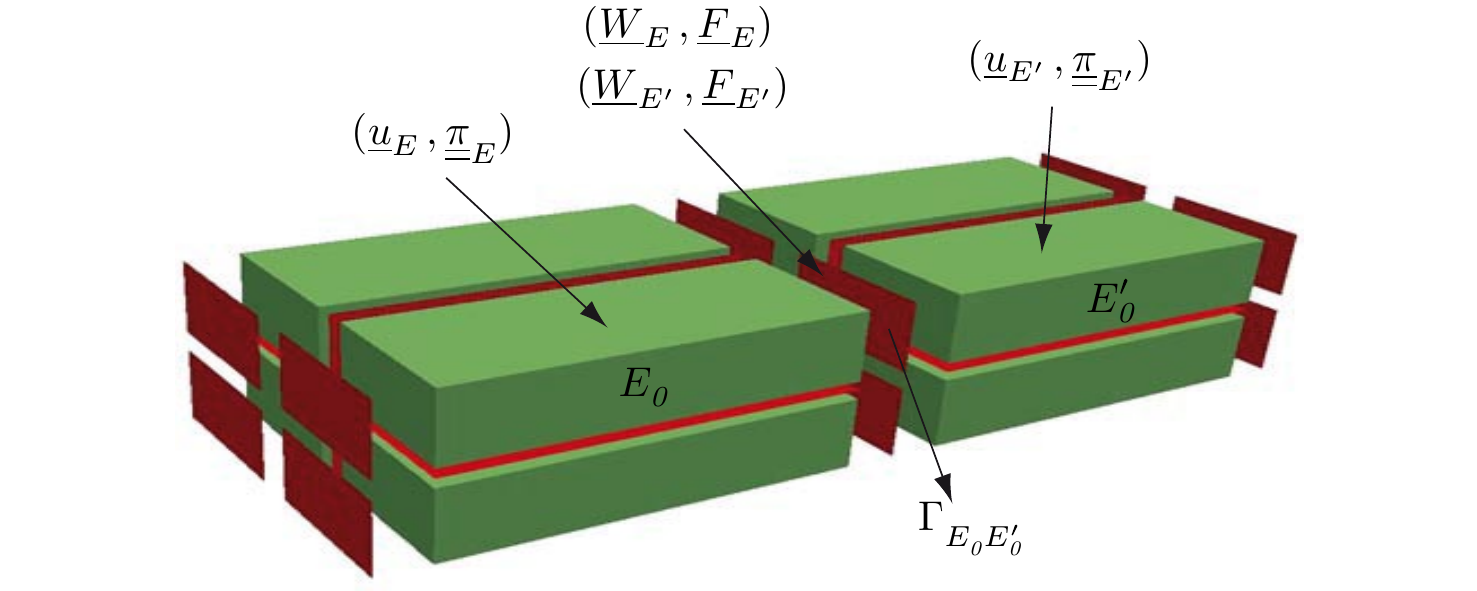}
       \caption{Unknown fields on the interfaces and substructures}
       \label{fig:sst}
\end{figure}

The laminated composite structure $\struct$ is divided into substructures and interfaces as shown in Figure \ref{fig:decomp_sst_interfaces}. Each of these mechanical entities possesses its own kinematic and static unknown fields related by its constitutive law. The substructuring process is governed by the objective of matching the domain decomposition interfaces with the cohesive material interfaces, so that each substructure belongs to a unique ply $\ply$ and its behavior is geometrically nonlinear. In the initial configuration $\Confio$, a substructure $\E$ defined in domain $\OmegaEo$ is connected to an adjacent substructure $\Ep$ through an interface $\GammaEEpo=\partial \OmegaEo \cap \partial \OmegaEpo$ (see Figure \ref{fig:sst}). The surface entity $\GammaEEpo$ applies force distributions $\FEo$, $\FEpo$ and displacement distributions $\WE$, $\WEp$ to $\E$ and $\Ep$ respectively. Let $\GammaEo = \bigcup_{\Ep \in \struct} \GammaEEpo$.

For a substructure $\E$ such that $\GammaEo \cap ( \dOmegaFdo \cup \dOmegaUdo) \neq \emptyset$, the boundary condition $(\Fdo,\Udo)$ is applied through a boundary interface $\GammaEdo$.

Let $\uE$ be the displacement field, $\GLE$ the Green-Lagrange strain tensor, $\GLpE$ the Lagrangian strain rate, $\piE$ the second Piola-Kirchhoff stress tensor, $\FFE$ the deformation gradient tensor and $J_{E} = det(\FFE)$ the Jacobian of the motion at each point of substructure $\E$.  At each point of interface $\GammaEEpo$, the displacement field is defined as $\WE$, the displacement gap as $[\WE] =\WEp-\WE$, the mean value as $\WmoyEE = \frac{1}{2}(\WEp+\WE)$ and the deformation gradient tensor as $\FFEE = \grado \WmoyEE+ \underline{\underline{I}}_{d}$.

Then, the substructured quasi-static problem consists, at each step of the time integration, in finding $s=
(s_\E)_{\E \in \mathbf{E}} $, where $s_\E = (\uE, \piE, \WE, \FEo )$ verifies the following equations:
\begin{itemize}
\item Mass conservation of substructure $\E$:
\begin{equation} \label{eq:mass_E}
\rho_E \, J_E = \rho_{E_\textit0} \;,\; \text{over} \; \OmegaEo
\end{equation}
\item Nonlinear kinematic admissibility of substructure $\E$:
\begin{equation} \label{eq:kine_add_E1}
\GLE= \frac{1}{2}\left(\grado \uE + \gradot \uE + \gradot \uE \grado \uE\right)
\ ,\ \text{over} \ \OmegaEo
\end{equation}
\begin{equation} \label{eq:kine_add_E2}
\uE_{| \dOmegaEo} = \WE_{|\GammaEo} \;,\; \text{over} \; \GammaEEpo
\end{equation}
\item Global nonlinear equilibrium of substructure $\E$:
\begin{multline}
\forall (\uEstar,\WEstar) \in \suEO \times \sWEO \;, \text{such that} \; \uEstar_{| \dOmegaEo} = \WEstar_{| \GammaEo},
\\
\int_\OmegaEo  \piE: \EpsilonpEstar \ d\Omegao = \int_\OmegaEo \rho_{\E} \; \fd \cdot \uEstar \, d \Omegao + \int_{\GammaEo} \FEo \cdot \WEstar \, d \Gammao  \;,
\label{eq:eq_sst}
\end{multline}
where $\EpsilonpEstar =\frac{1}{2}(\grado \uEstar + \gradot \uEstar + \gradot \uE \grado \uEstar + \gradot \uEstar \grado \uE)$.
\item Hyperelastic orthotropic behavior of substructure $\E$:
\begin{equation} \label{eq:rdc_E}
\piE = \frac{\partial \psi}{\partial \GLE} \;,\; \text{over} \; \OmegaEo\;,\;
\end{equation}
where $\psi$ is the stored energy function or elastic potential per unit of undeformed volume. In this first study, we use $\psi=\frac{1}{2} \KEo \, \GLE: \GLE$.
\item Constitutive equation of interface $\GammaEEpo$:
\begin{equation} \label{eq:rdc_inter}
\displaystyle {\mathcal{R}}_{\E\Ep} (   [\WE] \, , \, \FEo \, , \, \FEpo\, ,\, \FFEE) = \underline0 \;,\; \text{over} \; \GammaEEpo\in \GammaEo \;.
\end{equation}
\item Behavior of the interface at the boundary $\GammaEdo$:
\begin{equation} \label{eq:rdc_inter_boundary}
\mathcal{R}_{E_{d_\textit0}}( \WE, \FEo) = \underline0 \;,\;  \text{over} \; \GammaEdo\;.
\end{equation}
\end{itemize}
The formal relation $\mathcal{R}_{\E\Ep}=\underline0$, called the ``interface behavior'' and defined over $\GammaEEpo$, can be made explicit in the three cases we are concerned with:
\begin{itemize}
\item Perfect interface:
\begin{equation}
\left\{ \begin{array}{l}
\FEo + \FEpo = \underline0 \\
\text{} [\WE]  = \underline0
\end{array} \right.
\label{eq:rdc_inter_perfect}
\end{equation}
\item Cohesive interface:
\begin{equation}
\left\{ \begin{array}{l}
\FEo + \FEpo = \underline0
\\
\displaystyle \mathcal{A}_{PP'}([\WE] \, ,\, \FEo \, ,\, \FFEE ) = \underline0
\end{array} \right.
\label{eq:rdc_inter_cohesive}
\end{equation}
\item Unilateral contact interface (without friction):
\begin{equation}
\left\{ \begin{array}{l}
\FEo + \FEpo = \underline0 \\
\normalply \cdot [\WE]   \geqslant 0 \;\; \text{and} \;\; \normalply \cdot \FEo \geqslant 0\\
(\normalply \cdot [\WE] )(\normalply \cdot \FEo) = 0 \\
\mathbf{P} \FEo = \mathbf{P} \FEpo = \underline0 \\
\end{array} \right.
\label{eq:rdc_inter_contact}
\end{equation}

where substructures $\E$ and $\Ep$ belong respectively to plies $\plyo$ and $\plypo$; operator $ \mathcal{A}_{PP'}$ was introduced in Section \ref{section:interface_law}; $\normalply$ was defined in Eq. \eqref{eq:deformed_normal} and $\mathbf{P}$ is the corresponding tangential projection operator.

As already mentioned at the end of Section \ref{section:interface_law}, contact in the delaminated interfaces is detected only between points which had the same initial position.

\textit{Remark:} Usually, the admissibility of forces is written in the deformed configuration; here, it can be easily converted to the initial configuration thanks to the geometric coincidence of the deformed interface (Assumption \ref{hyp_inter}).
\end{itemize}

\section{The numerical resolution strategy}
\label{section:strategy}

\subsection{The macroscopic scale}
\label{subsubsec:macro}

In order to ensure the scalability of the method, one can solve a coarse global linear problem. The definition of the macroscopic fields required to formulate this problem refers only to the interface's unknowns. The action/reaction principle $\FEo + \FEpo = \underline{0}$ is verified regardless of the interface's behavior. The objective of the macroscopic problem is to ensure that part of this equation is verified at any time:
\begin{equation}\label{eq:macro_eq}
 \int_{\GammaEo}  (\FEo + \FEpo)\cdot\WMstar\, d \Gammao= 0,\ \forall \WMstar \in \sWM\;,
\end{equation}
where the displacement macrospace $\sWM$ and its dual space $\sFM$ are parameters of the method. These subspaces are common to neighboring substructures and induce a separation of the interface quantities which is made unique by the uncoupling of the virtual works:
\begin{multline}
\forall  (\FEo,\WE)  \in \sFE\times\sWE, \quad \FEo = \FEoM + \FEom   \;,  \; \WE = \WEM + \WEm \;, \\
 \int\limits_{\GammaEo} \FEo\cdot\WE \, d \Gammao
=  \int\limits_{\GammaEo} \FEoM\cdot\WEM \, d \Gammao  + \int\limits_{\GammaEo} \FEom\cdot\WEm \, d \Gammao \;,
\end{multline}

\begin{figure}[t]
       \centering
       \includegraphics[width=0.65 \linewidth]{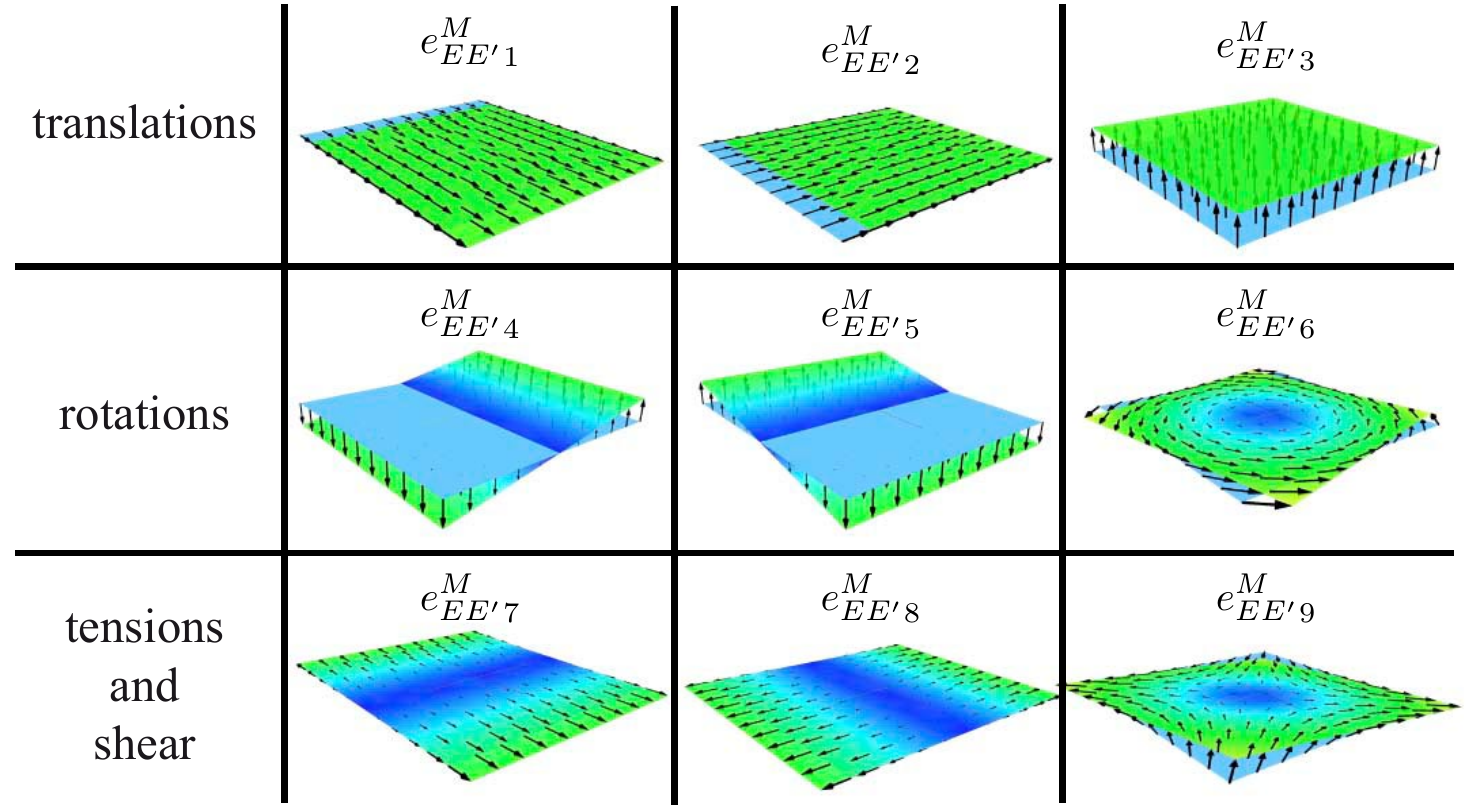}
       \caption{The linear macrobasis for a plane interface}
       \label{fig:base_macro}
\end{figure}

Usually, one chooses a common basis for the kinematic and static macroscopic fields of the interface. Numerical tests have shown that in order to ensure the numerical scalability of the method the macroscopic basis should extract at least the linear part of the interface forces (see Figure \ref{fig:base_macro}). Indeed, this macroscopic space contains the part of the interface fields with the longest wavelength. Consequently, according to Saint-Venant's principle, the micro complement has only local influence.

\subsection{The iterative algorithm}\label{subsec:resol}

In this section, the iterative LATIN algorithm, which enables one to deal with nonlinear problems, is adapted to the resolution of the geometrically nonlinear substructured reference problem with nonlinearities localized at the interfaces. The finite element method is used to discretize the equations.

The equations of the problem are divided into two groups:
\begin{itemize}
    \item Admissibility of the substructures and macroscopic admissibility of the interfaces:
    \begin{enumerate}
         \item[-] mass conservation of the substructures, Eq. \eqref{eq:mass_E};
        \item[-] nonlinear kinematic admissibility of the substructures, Eq. (\ref{eq:kine_add_E1}, \ref{eq:kine_add_E2});
        \item[-] nonlinear static admissibility of the substructures, Eq. \eqref{eq:eq_sst};
        \item[-] behavior of the substructures, Eq. \eqref{eq:rdc_E};
        \item[-] macroscopic admissibility of the interfaces (after the linearization of the previous equations), Eq. (\ref{eq:macro_eq}).
    \end{enumerate}
    \item Local (nonlinear) equations at the interfaces:
    \begin{itemize}
        \item[-] interface behavior, Eq. (\ref{eq:rdc_inter}, \ref{eq:rdc_inter_boundary}).
    \end{itemize}
\end{itemize}

\begin{figure}
       \centering
       \includegraphics[width=.65 \linewidth]{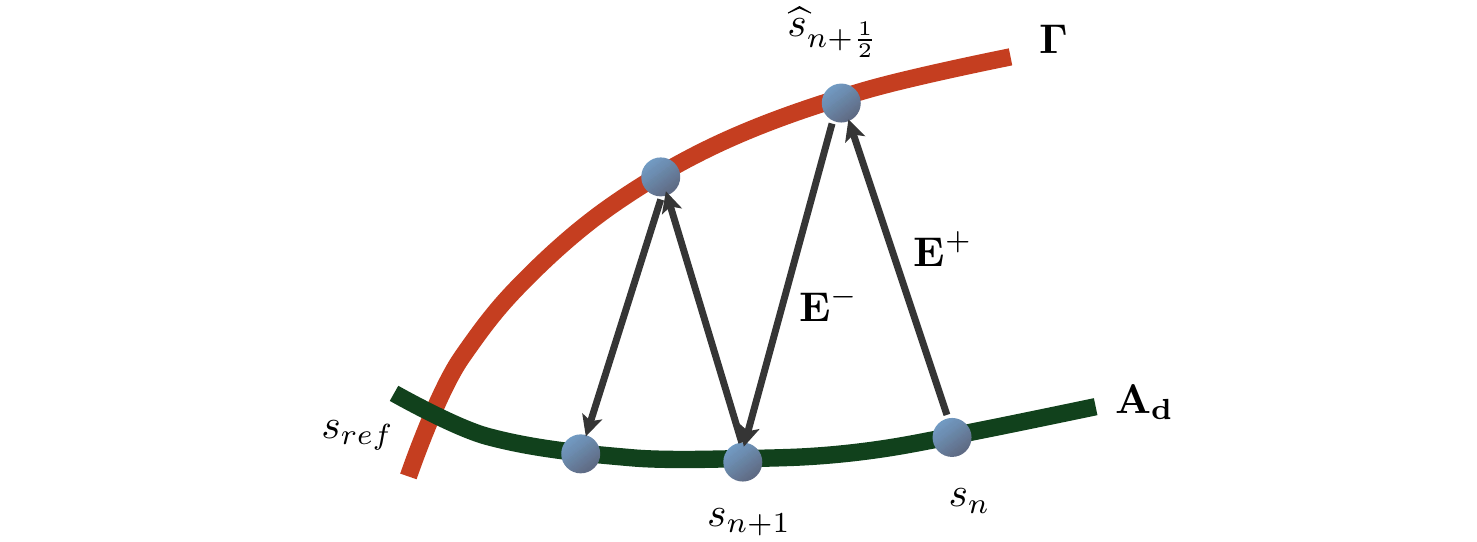}
       \caption{Schematic representation of the LATIN iterative algorithm}
       \label{fig:algo_latin}
\end{figure}

The interface solutions $s = (s_\E)_{\E \in \struct} = (\WE , \FEo )_{\E \in \struct}$ of the first set of equations belong to space $\sAd$, while the interface solutions $\widehat{s} = (\widehat{s}_\E)_{\E \in \struct} = (\WchapE , \FchapEo )_{\E \in \struct}$ of the second set of equations belong to $\sGamma$. The converged interface solution $s_{ref}$ is such that:
\begin{equation}
s_{ref} \in \sAd \cap \sGamma \;.
\end{equation}
The resolution process consists in seeking the interface solution $s_{ref}$ alternatively in these two spaces: first, one finds a solution $s_n$ in $\sAd$, then a solution $\widehat{s}_{n+\frac{1}{2}}$ in $\sGamma$. In order for the two problems to be well-posed, one introduces the search directions $\sEp$ and $\sEm$ which link the solutions $s$ and $\widehat{s}$ during the iterative process (see Figure \ref{fig:algo_latin}).
Thus, an iteration of the LATIN algorithm consists of two nonlinear stages which are described in detail below.

\textit{Remark:} The LATIN approach \cite{Ladeveze99,Ladeveze03a} is a general computational strategy for the resolution of time-dependent nonlinear problems which operates over the entire space-time domain. In our case, time is irrelevant and the capabilities of the LATIN method are not fully exploited.
The LATIN strategy is also based on the idea of separating the difficulties and dealing with global linear equations and local nonlinear equations independently. In this work, we consider a global nonlinear stage (called the admissibility stage) which is solved as a series of global linear equations.

\subsubsection{The local stage:}
In the local stage, the following local problems are solved at each point of interfaces $\GammaEEpo$:
\begin{equation}
\label{eq:local_problem}
\begin{array}{l}
\text{Find} \ (\FchapEo,\WchapE,\FchapEpo,\WchapEp) \text{ such that: } \
\left\{ \begin{array}{l}
\displaystyle \mathcal{R}_{\E\Ep}([\WchapE],\FchapEo,\FchapEpo, \FFchap_{E^{\phantom{'}}_\textit0E'_\textit0}) = \underline0  \\
\displaystyle  (\FchapEo-\FEo) -\kpEo (\WchapE-\WE) = \underline0 \\
\displaystyle  (\FchapEpo-\FEpo) -\kpEpo (\WchapEp-\WEp) = \underline0
\end{array} \right.
\end{array}
\end{equation}
The last two equations of this system define the search direction $\sEp$. In the case of a cohesive interface, Problem \eqref{eq:local_problem} is nonlinear and is solved using a modified Newton-Raphson algorithm.

The following local linear problems are also solved at each point of boundary interfaces $\GammaEdo$:
\begin{equation}
\begin{array}{l}
\text{Find} \ (\FchapEo,\WchapE) \ \text{such that: }
\left\{ \begin{array}{l}
\displaystyle \mathcal{R}_{E_{d_\textit0}}(\WchapE,\FchapEo) = \underline0 \\
\displaystyle  (\FchapEo-\FEo) -\kpEo (\WchapE-\WE) = \underline0
\end{array} \right.
\end{array}
\end{equation}

\subsubsection{The admissibility stage:}

The admissibility stage consists in solving the nonlinear equations in each substructure under the constraint of macroscopic admissibility of the interfaces. This stage is carried out using an iterative Newton-Raphson procedure; therefore, at each iteration, the macroscopic admissibility is prescribed on the linearized system of equations in the substructures.

In absence of macroscopic admissibility, the search direction $\mathbf{E}^-$ which couples the interface displacement and the force fields from the linear stage is:
\begin{equation}
\label{eq:ddr_global_mono}
(\FEo-\FchapEo) +\kmEo (\WE - \WchapE) = \underline0 \;,\; \text{over} \; \GammaEEpo\;.
\end{equation}
The monoscale version of the nonlinear problem in each substructure becomes:
\begin{multline}
\label{eq:pb_sst_mono_NL}
\forall (\uEstar,\WEstar) \in \suEO \times\sWEO,
\\
\underbrace{ \int_\OmegaEo \KEo \, \GLE (\uE): \EpsilonpEstar \, d\Omegao + \int_\GammaEo \kmEo \, \WE \cdot \WEstar \, d \Gammao }_{-P_{int}(\uE,\WE )}
\\
=  \underbrace{ \int_\OmegaEo \rho_\E \fd \cdot \uEstar \, d \Omega + \int_\GammaEo (\FchapEo + \kmEo \, \WchapE) \cdot \WEstar \, d \Gammao }_{P_{ext}}\;.
\end{multline}

This problem is solved using a Newton-Raphson algorithm. At each iteration $i$, the correction term $({}^i\delta\uE, {}^i\delta\WE)$ is calculated by linearizing \eqref{eq:pb_sst_mono_NL}:
\begin{multline}
\label{eq:pb_sst_mono_L}
\forall (\uEstar,\WEstar) \in \suEO \times\sWEO,
\\
 \int_\OmegaEo \KEo \, \GLE ({}^i\uE): \left(\gradot {}^i\delta\uE \; \grado \uEstar + \gradot \uEstar \; \grado {}^i\delta\uE \right) \, d\Omegao
 \\
 +  \int_\OmegaEo \KEo \, \EpsilonpE({}^i\delta\uE): \EpsilonpE({}^i\uEstar) \ d\Omegao + \int_\GammaEo \kmEo \, {}^i\delta\WE \cdot \WEstar \, d \Gammao =  P_{ext} +P_{int}({}^i\uE, {}^i\WE) \;,
\end{multline}
where $\EpsilonpE(\underline{v}) =\frac{1}{2}(\grado \underline{v} + \gradot \underline{v} + \gradot {}^i\uE \grado \underline{v}  + \gradot \underline{v} \grado {}^i\uE )$, $({}^i\uE, {}^i\WE)$ is known from the previous iteration, and $({}^{i+1}\uE,{}^{i+1}\WE) = ({}^i\uE, {}^i\WE) + ({}^i\delta\uE, {}^i\delta\WE)$.

Then, at each iteration, macroscopic admissibility is applied to linearized problem \eqref{eq:pb_sst_mono_L}. The satisfaction of the macroscopic equilibrium of interface forces \eqref{eq:macro_eq} suffices to ensure the scalability of the method.

Condition \eqref{eq:macro_eq} is incompatible with search direction \eqref{eq:ddr_global_mono}; hence this search direction is weakened and verified ``as best can be'' under the macroscopic constraint. Technically, this is achieved by using a Lagrangian multiplier whose stationarity leads to a modified local search direction:
\begin{multline} \label{eq:ddr_loc}
\forall \WEstar \in \sWEO, \quad \int_\GammaEo  ({}^{i+1}\FEo-   \FchapEo) \cdot \WEstar \, d\Gammao
\\
+ \int_{\GammaEo}  \left( \kmEo \ ({}^{i+1}\WE-\WchapE)   -   \kmEo \, \Wtilde \right) \cdot \WEstar \, d\Gammao= 0 \;.
\end{multline}

In order to solve this system, a relation linking ${}^{i+1}\FEoM$ and ${}^{i+1}\WE$ is derived from the subdomain equilibrium plus the modified search direction (\ref{eq:pb_sst_mono_L}, \ref{eq:ddr_loc}). The linear problem to be solved in order to find $({}^i\delta\uE, {}^i\delta\WE)$ becomes:
\begin{multline}
\label{eq:pb_micro}
\forall (\uEstar,\WEstar) \in \suEO\times\sWEO,
\\
 \int_\OmegaEo \KEo \, \GLE ({}^i\uE): \left(\gradot {}^i\delta\uE \; \grado \uEstar + \gradot \uEstar \; \grado {}^i\delta\uE \right) \, d\Omegao
 \\
 +  \int_\OmegaEo \KEo \, \EpsilonpE({}^i\delta\uE): \EpsilonpE(\uEstar) \, d\Omegao + \int_\GammaEo \kmEo \, {}^i\delta\WE \cdot \WEstar \, d \Gammao =  {}^i\tilde{P}_{ext} +P_{int}({}^i\uE, {}^i\WE) \;,
\end{multline}
where $\displaystyle {}^i\tilde{P}_{ext} = P_{ext} +  \int_\GammaEo \kmEo \, \Wtilde \cdot \WEstar \, d \Gammao$. One can prove that if $\KEo$ and $\kmEo$ are symmetric, positive definite operators, then Eq. \eqref{eq:pb_micro} is well-defined and has a unique solution. Due to the linearity of Eq. \eqref{eq:pb_micro}, one can define a linear relation between the interface displacements and the loading:
\begin{multline}
\forall  \  \FEstar \in  \sFEO  , \quad \int_\GammaEo  {}^{i+1}\WE \cdot \FEstar \, d\Gammao \\
= \int_\GammaEo \left(\HE (\FchapchapEo + \kmEo \WtildeE ) + \WEdc + {}^{i}\WE \right) \cdot \FEstar d\Gammao\;,
\end{multline}
where $\displaystyle \FchapchapEo=\FchapEo + \kmEo \, \WchapE$, and operator $\HE$ is the dual Schur complement of substructure $\E$ modified by the search direction at iteration $i$ (which depends on the geometric configuration of the previous iteration $i$), while $\WEdc$ results from the condensation of the volume loading and of $P_{int}({}^i\uE,{}^i\WE)$ onto interface $\GammaEo$. ${}^{i}\WE$ are the interface displacements of the iteration $i$.

The corresponding interface forces are obtained using the modified search direction \eqref{eq:ddr_loc} and projected onto the macroscopic space:
\begin{multline}
\label{eq:comp_homo}
\forall \  \WMstar \in  \sWM , \\
\int_\GammaEo {}^{i+1}\FEo \cdot \WMstar \, d\Gammao
= \int_\GammaEo ( \LEM \ \Wtilde + \FtildeEo)  \cdot   \WMstar \, d\Gammao \, , \,
\end{multline}
where:
\begin{multline*}
\forall \displaystyle \WMstar \in  \sWM, \\
\int_\GammaEo \LEM \Wtilde  \cdot  \WMstar d \Gammao = \int_\GammaEo ( \kmEo - \kmEo \HE \kmEo ) \Wtilde \cdot \WMstar d\Gammao
\\
\int_\GammaEo \FtildeEo \cdot  \WMstar d \Gammao = \int_\GammaEo  (\FchapchapEo
- \kmEo ( \HE \FchapchapEo + \WEdc + {}^{i}\WE ))  \cdot  \WMstar d\Gammao \,.
\end{multline*}
$\LEM$ is classically viewed as the homogenized behavior of substructure $\E$ and is calculated explicitly for each substructure by solving local subproblems \eqref{eq:pb_micro}, taking the vectors of the macroscopic basis as boundary conditions over $\GammaEo$. Clearly, $\LEM$ depends on the geometric configuration of the previous iteration.

Finally, Relation \eqref{eq:comp_homo} is introduced into Eq. \eqref{eq:macro_eq}, which expresses the admissibility of the macroforces, leading to what is called the macroscopic problem:
\begin{multline}
\label{eq:pb_macro}
\forall \ \Wtildestar \in \sWMad, \quad
\sum_{\E \in \struct} \int_\GammaEo \LEM \, \Wtilde \cdot \Wtildestar \, d\Gammao
\\ =  \sum_{\E \in \struct} \int_\dOmegaFdo \Fdo \cdot \Wtildestar \, d\Gammao - \sum_{\E \in \struct} \int_\GammaEo  \FtildeEo \cdot \Wtildestar \, d\Gammao
\end{multline}
The macroscopic problem is discrete by nature. Therefore, it has an algebraic form $\displaystyle {}^i\mathbf{L^M} \ \Wtilde = {}^{i+1}\underline F^M$, where $\Wtilde$ is the vector of the components of the Lagrange multiplier in the macroscopic basis.

The right-hand side of Eq. \eqref{eq:pb_macro} can be interpreted as a macroscopic static residual obtained from the calculation of a monoscale linear stage. In order to derive this term, Problem \eqref{eq:pb_micro} must be solved independently in each substructure, setting $\Wtilde$ to zero.
The resolution of the macroscopic problem of Eq. \eqref{eq:pb_macro} leads to the global knowledge of Lagrange multiplier $\Wtilde$, which is finally used as a prescribed displacement to solve the substructure-independent problems of Eq. \eqref{eq:pb_micro}.

The resolution of \eqref{eq:pb_micro} in the substructures is carried out using the finite element method.

\vspace{5mm}
\noindent \textit{Remarks:}
\begin{itemize}
\item The convergence of the algorithm can be improved by introducing a relaxation stage after the admissibility stage. The admissibility solution $s_n$ is renamed $\breve{s}_n$; then the relaxed solution $s_n$ is defined by:
\begin{equation}
s_n = \mu \breve{s}_n + (1-\mu)s_{n-1}\;,
\end{equation}
where $\mu$ is a relaxation parameter usually taken equal to $0.8$.
\item The LATIN error indicator adopted here was successfully used in \cite{Kerfriden09,Allix10}. This criterion is based on a measure of the non-satisfaction of the constitutive laws of Eq. (\ref{eq:rdc_inter}, \ref{eq:rdc_inter_boundary}) in the admissibility stage, since these are the only equations which are not verified at this stage. More precisely, each time a solution $s_n \in \sAd$ is obtained, an indicator of the convergence of the algorithm is calculated by integrating the local residuals of the interface behavior over the structure.
\end{itemize}

\section{Analysis of the parameters of the algorithm in the case of slender structures}
\label{sec:analysis}

The search direction parameters of the local stage $(\kpEo)_{(\E \in \mathbf{E})}$ and of the admissibility stage $(\kmEo)_{(\E \in \mathbf{E})}$ are symmetric, positive definite operators which represent the influence of the neighboring subdomains and interfaces, such as interface impedances or Schur complements of the rest of the structure $\Omegao \setminus \OmegaEo$. It was empirically shown in previous studies \cite{Champaney99, Ladeveze01} that there is an optimum set of these operators which depends on the interface behavior.

For the monoscale approach (which does not include the satisfaction of the macroscopic problem) applied to massive isotropic and homogenous structures with contact or perfect interfaces, the basic setting is the scalar approximation $\kpEo=\kmEo \simeq E/L_\textit0$, where $E$ is Young's modulus and $L_\textit0$ is a characteristic length of the structure. The multiscale approach enables most of the effect of the search direction to be localized, so the classical setting becomes $\kpEo=\kmEo \simeq E/L_\GammaEo$, where $L_\GammaEo$ is a characteristic length of the interface.

In the case of cohesive interfaces with damage, the local search directions $\kpEo$ of the cohesive interfaces must be set to infinite values, because values which are too small can lead to the stagnation or the divergence of the algorithm, as explained in \cite{Kerfriden09}. This choice also enables the interface's quantities to be calculated directly in the local stage, so local Newton iterations are unnecessary.
Concerning search direction $\kmEo$, the optimum value would be twice the actual interface stiffness $2 k^0 (1-d)$ ($k^0$ denoting the undamaged interface's local stiffness), which would be equivalent to prescribing the exact interface behavior as an interface condition in the admissibility stage {\cite{Kerfriden09}}. Unfortunately, the use of this value would require the operators to be updated very frequently, which would be expensive; instead, a monitoring strategy has been proposed in \cite{Kerfriden09}.

Our particular study, which focuses on multiple buckling-delamination interactions in composite laminates, is characterized by the occurrence of bending in very slender geometries because the plies are very thin ($\approx 0.125$~mm) and small delaminated areas can have high slenderness coefficients ($L_{delaminated}/h_{ply} \gg 100$). As a result, classical values for massive structures are completely inadequate and the following difficulties need to be tackled:
\begin{itemize}
\item Concerning the search direction for perfect interfaces, the scalar $E/L_\GammaEo$ nearly equals the stiffness of the neighbors (disregarding the geometry of the structure), which could affect performance in the case of slender structures, especially in bending and in buckling ({see Section \ref{section:bending}}).
\item The occurrence of buckling leads to major changes in the deformed configuration which require additional continuity conditions in order to ensure convergence, {as proven in Section \ref{section:buckling}.}
\item In the case of multiple buckling-delamination interactions, it is possible for delaminated surfaces to separate, regain contact or remain in contact throughout the evolution, so during the calculations the same surface can go through different states which are unpredictable.
Unfortunately, in the case of slender structures and subdomains, there may be different suitable search directions depending on the state (open or closed), and an inadequate value could lead to nonphysical solutions or to stagnation of the iterative process. {See Section \ref{section:ddr_contact} for the proposed remedies.}
\item {The use of optimum search directions $\kmEo$ for cohesive interfaces \cite{Kerfriden09} requires the operators to be updated and reconstructed according to the damage state at each point of the cohesive interface, which leads to an increase in CPU time. In the procedure proposed in \cite{Kerfriden09}, the search directions are updated only when an interface becomes fully damaged. Unfortunately, the efficiency of this procedure deteriorates when the size of the substructures increases, and becomes even worse in the case of slender subdomains with large cohesive interfaces (because the updating is performed only once the entire interface is damaged, which can take several time steps). See Section \ref{section:ddr_cohesive_interfaces} for the recommended solution.}
\end{itemize}

Finding the most favorable values requires not only additional identifications in order to take into account the geometry of the structure, but also the introduction of a new step: the micro/macro separation of the search direction.

If one divides search direction $\mathbf{E}^-$ of the admissibility stage into a macroscopic part ${\mathbf{E}^-}^\mathbf{M}$ and a microscopic part ${\mathbf{E}^-}^\mathbf{m}$, Eq.~\eqref{eq:ddr_loc} can be rewritten as:
\begin{multline}
\label{eq:ddr_macro_sep}
\forall \ \WMstar  \in \sWM, \quad \\
\int_\GammaEo ( \FEo-\FchapEo )\cdot \WMstar \, d\Gammao + \int_\GammaEo \kmEoM \ (\WE-\WchapE -\underline{\widetilde{W}}_{E_\textit0}^M )\cdot \WMstar \, d\Gammao= 0
\end{multline}
\begin{multline}
\label{eq:ddr_micro_sep}
\forall \ \Wmstar \in \sWm, \quad
\int_\GammaEo ( \FEo-\FchapEo ) \cdot \Wmstar \, d\Gammao + \int_\GammaEo  \kmEom \ (\WE-\WchapE) \cdot \Wmstar \, d\Gammao= 0
\end{multline}
where $\sWm$ is the space orthogonal to $\sWM$ with respect to the inner product $L^2(\GammaEo)$.

Parameter $\kmEoM$ represents the stiffness of the interface for the macroscopic problem. Its optimum value, if it exists, is the homogenized interface behavior: for perfect interfaces (infinite stiffness), this parameter must be as large as possible; for homogenous elastic interfaces (\emph{e.g.}  cohesive interfaces with constant damage), it is related to the current stiffness $\kmEoM = 2k^0 (1-d)$ (where $k^0$ is the undamaged interface's local stiffness). Parameter $\kmEom$ is the micro part of the search direction and, classically, $\kmEom= E/L_\GammaEo$ ($\kmEom= 2k^0 (1-d)$) is considered to be a good starting value for perfect interfaces (elastic interfaces).

In order to deal with the difficulties one at a time, we will seek the optimum parameters independently for the different types of interface behavior using academic examples of 3D slender structures. In the following sections, separate analyses will be presented for perfect interfaces in bending and buckling (Section \ref{section:ddr_perfect_interfaces}), contact interfaces (Section \ref{section:ddr_contact}) and cohesive interfaces (Section \ref{section:ddr_cohesive_interfaces}).

All the analyses were carried out using a fully parallel C++ program, taking advantage of the three scales of the mixed domain decomposition method proposed. The transfers of data among the different processors required for the parallel computations were performed using the MPI libraries. Each processor was assigned to a set of connected substructures and their interfaces and used to calculate the associated operators and solve the local problems. This was achieved technically thanks to a METIS routine and helped reduce the number of interfaces among the processors.

\subsection{Perfect interfaces}
\label{section:ddr_perfect_interfaces}

Here, we present a study of the search directions for slender structures. The analysis concerns a cantilever plate in bending under the assumption of small perturbations and shows how a thin geometry affects the convergence rate and the scalability of the strategy (Section \ref{section:bending}). The proposed improvements as a result of the bending study will be tested in the case of a buckling example in Section \ref{section:buckling}.

\subsubsection{Bending:}
\label{section:bending}
\begin{figure}[ht]
       \centering
       \includegraphics[width=0.65 \linewidth]{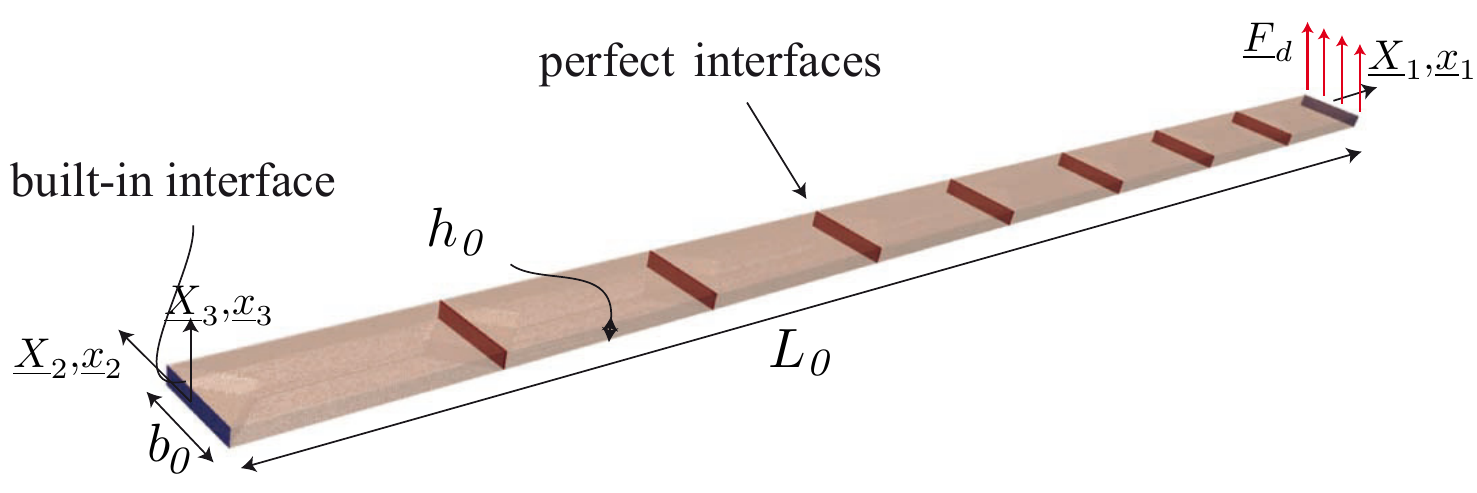}
       \caption{The 3D cantilever plate in bending: geometry, substructures and interfaces}
       \label{fig:flexion_confi}
\end{figure}

{The problem being considered is that of a slender plate which is built-in along one side and subjected to a bending load $\underline F_d$ in the form of a surface force distribution along the opposite side, as shown in Figure \ref{fig:flexion_confi}. The material is isotropic and homogeneous, the geometry was fixed to: $L_\textit0=20$~mm, $h_\textit0=0.2$~mm and $b_\textit0=1$~mm. The substructuration was made only along the $X_1$-direction, modifying the number of substructures to study the convergence rate.} Figure \ref{fig:flexion_confi} shows the substructuration in 8 subdomains. The whole mesh totaled 1.8 million DOFs with 10 linear wedge elements through the thickness. All the calculations were made using $\kpEo=\kmEom$.

Table \ref{tab:ddr_flexion} presents the number of iterations as a function of the number of substructures. Row A corresponds to the use of the classic values $\kmEoM = \kmEom = E/L_\GammaEo$, while Row B corresponds to continuous macrodisplacements $\kmEoM \to \infty$, keeping $\kmEom = E/L_\GammaEo$. In these two cases the convergence rate depends on the number of substructures and deteriorates when this number increases, even if the continuity of the macrodisplacements is enforced.

\begin{table}[h] \footnotesize
\begin{center}
    \begin{tabular}{ | c | c | c | c | c | c | }
    \hline
    & number of substructures & 8 & 16 & 32 & 64 \\
    \hline
    \multirow{3}{*}{A} & number of LATIN iterations until convergence ($0.1\%$) & \multirow{3}{*}{15} & \multirow{3}{*}{15} & \multirow{3}{*}{22} & \multirow{3}{*}{38}    \\
    & without continuity of the macrodisplacements & & & &   \\
    & and without anisotropic search directions &   & & & \\
    \hline
    \multirow{3}{*}{B} & number of LATIN iterations until convergence ($0.1\%$) & \multirow{3}{*}{14} & \multirow{3}{*}{15} & \multirow{3}{*}{21} & \multirow{3}{*}{36}    \\
    & with continuity of the macrodisplacements &   & & &\\
    & and without anisotropic search directions &   & & &\\
    \hline
    \multirow{3}{*}{C} & number of LATIN iterations until convergence ($0.1\%$) & \multirow{3}{*}{3} & \multirow{3}{*}{6} & \multirow{3}{*}{9} & \multirow{3}{*}{15}    \\
    & without continuity of the macrodisplacements &   & & &\\
    & and with anisotropic search directions &  & & & \\
    \hline
    \multirow{3}{*}{D} & number of LATIN iterations until convergence ($0.1\%$) & \multirow{3}{*}{4} & \multirow{3}{*}{4} & \multirow{3}{*}{4} & \multirow{3}{*}{5}    \\
    & with continuity of the macrodisplacements &  & & & \\
    & and with anisotropic search directions &  & & & \\
    \hline
    \end{tabular}\\
\end{center}
\caption{Influence of the number of substructures on the convergence rate ($L_\textit0/h_\textit0 = 100$)}
\label{tab:ddr_flexion}
\end{table}

We propose an enhancement to the microscopic directions based on the fact that the local stresses and displacements in the normal and transverse directions are very different. Plate theory leads to the following relations between the orders of magnitude of the normal and shear stresses and between the orders of magnitude of the normal and transverse displacements:
\begin{equation*}
O\left(\frac{\tau _{X_1X_3}}{h_\textit0}\right)=O\left(\frac{ \sigma _{X_1X_1}}{L_\textit0}\right)\quad , \quad O\left(\frac{u_{X_1}}{h_\textit0}\right) = O\left(\frac{u_{X_3}}{L_\textit0}\right)\;.
\end{equation*}
Therefore, since the search directions relate tractions to displacements, they should take anisotropic values:
\begin{equation}
\label{eq:anisotropic_ddr}
 \frac{(\kmEom)_n}{(\kmEom)_t} =  \left(\frac{L_\textit0}{h_\textit0}\right)^2\;,
\end{equation}
where the ratio $L_\textit0/h_\textit0$ is a macroscopic quantity representing the whole structure, the normal value of $(\kmEom)_n = E/L_\GammaEo$ remaining constant.

Using these improved search directions and setting $\kmEoM = \kmEom$, the number of iterations dropped by more than $60\%$ (Row C of Table \ref{tab:ddr_flexion}), and the strategy became even more efficient and scalable when the continuity of the macrodisplacements $\kmEoM \to \infty$ was enforced (Row D of Table \ref{tab:ddr_flexion}).

\subsubsection{Buckling:}
\label{section:buckling}

\begin{figure}[!b]
       \centering
       \includegraphics[width=.65 \linewidth]{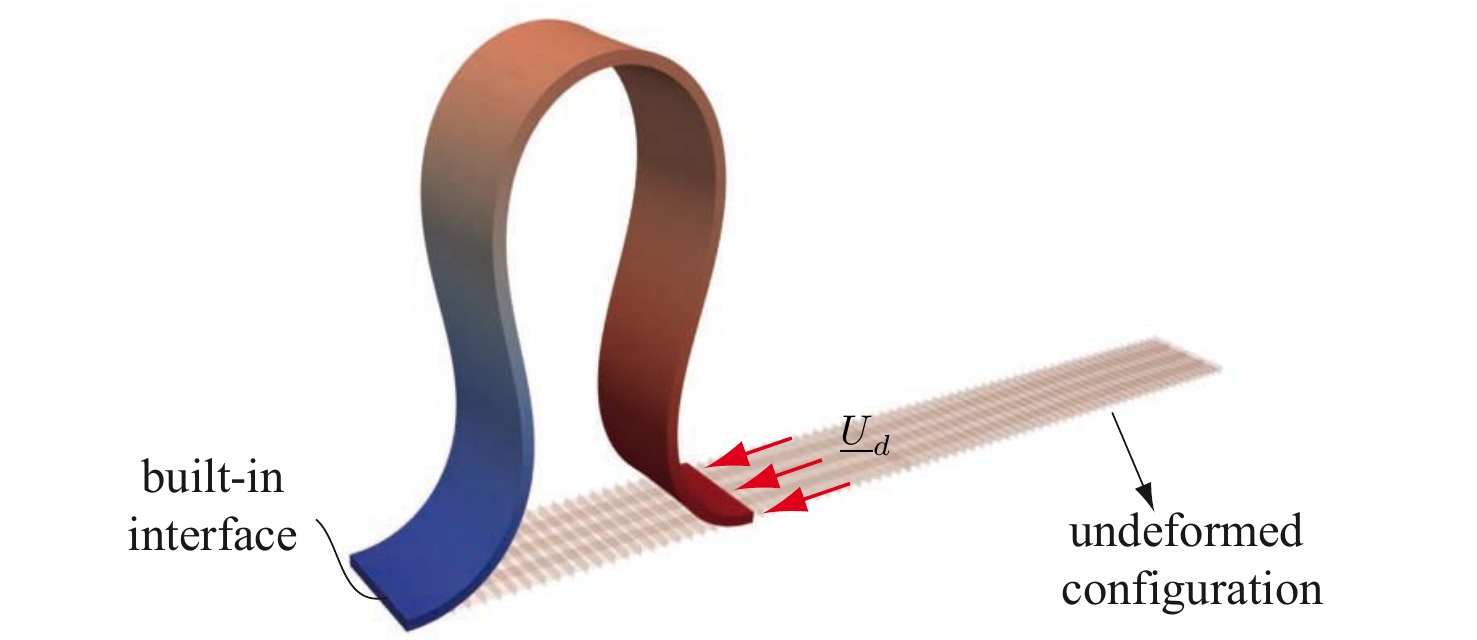}
       \caption{The initial configuration and the final deformed shape after the last time step}
       \label{fig:flam_deforme}
\end{figure}

In this section, we consider the geometrically nonlinear example of a plate built-in at both ends and subjected progressively to a negative end displacement, with a perturbation consisting of a central force. The data are: $L_\textit0 = 10$ mm, $h_\textit0 = 0.1$ mm, $b_\textit0 = 1$ mm, $E = 135,000$ MPa and $\nu = 0.3$. The geometry was divided into 640 substructures and 1,464 perfect interfaces (see Figure \ref{fig:flam_deforme}). The mesh totaled 2.2 million DOFs with 12 linear wedge elements through the thickness. The macroscopic problem and the supermacroscopic problem represented 13,176 DOFs and 372 DOFs respectively and were solved using 64 processors.

\begin{figure}[b]
       \centering
       \includegraphics[width=.65 \linewidth]{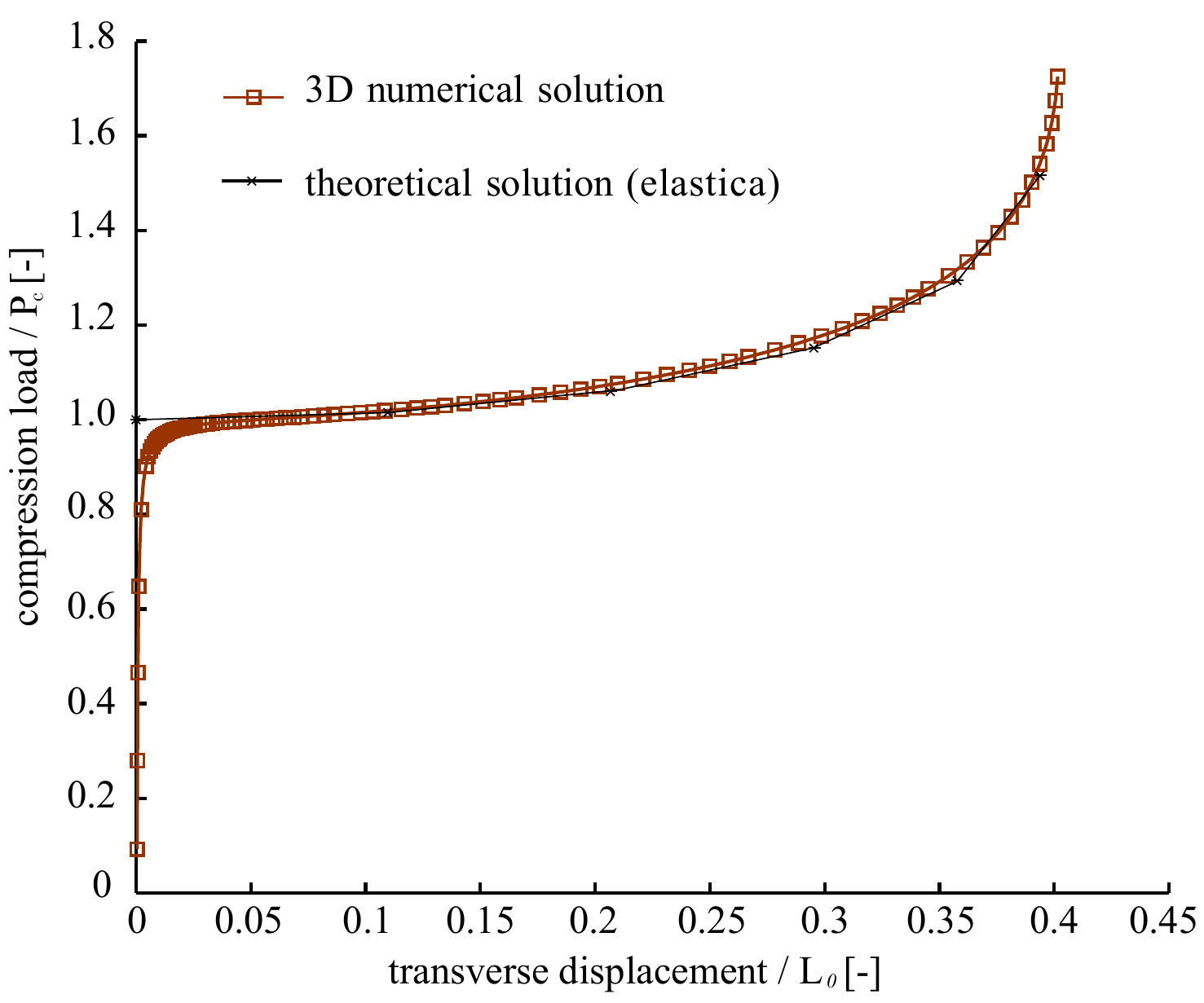}
       \caption{The load-displacement curve for a compressed cantilever plate}
       \label{fig:flam}
\end{figure}

The nonlinear buckling analysis was performed in 96 time steps, resulting in the evolution of the axial compression load as a function of the transverse displacement of the central perturbed point shown in Figure \ref{fig:flam}. This numerical response agrees perfectly with the theoretical response {given in \cite{Timo61}}, which has a critical Euler force equal to $P_c=\frac{4 \pi^2 E I}{L^2}=4.4$ N. Figure \ref{fig:flam_deforme} also shows the final deformed configuration. It is important to note that buckling began close to the $6^{th}$ time step.

\begin{figure}[b]
       \centering
       \includegraphics[width=.65 \linewidth]{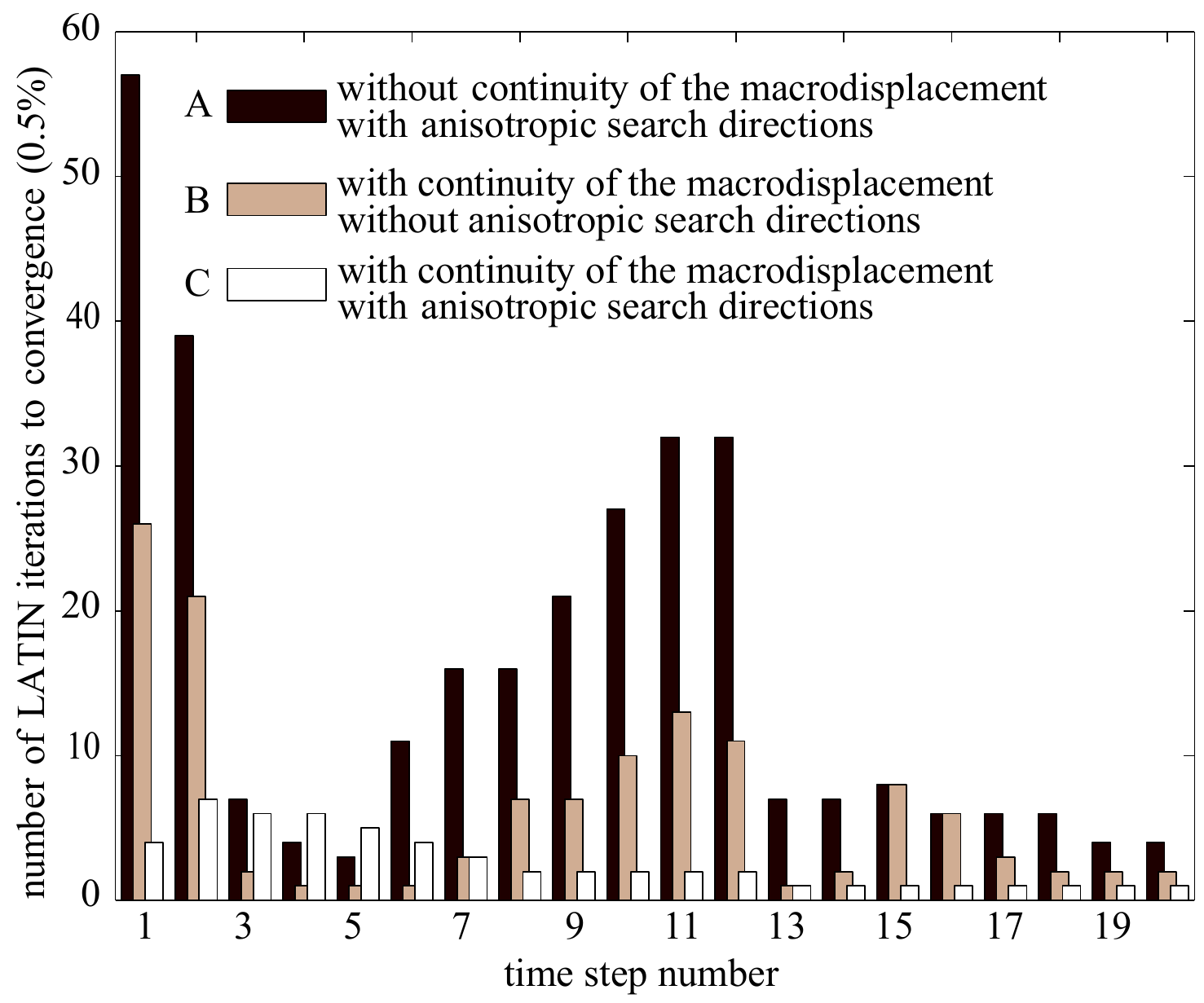}
       \caption{Influence of the search directions in the case of buckling}
       \label{fig:ddr_flam}
\end{figure}

{A search direction study similar to that presented before for bending under small perturbations was carried out for this buckling problem. Figure \ref{fig:ddr_flam} shows the numbers of LATIN iterations for each of the first 20 time steps corresponding to three different settings. In the first case, corresponding to the application of the values $\kmEoM = \kmEom = E/L_\GammaEo$, convergence was not reached after the $3^{rd}$ time step. Using the same values, but this time with anisotropic differentiation of the search direction, the method diverged before buckling (after the $6^{th}$ time step). Convergence could be achieved only by increasing the value $\kmEoM = \kmEom = E/L_\GammaEo$ at least a hundredfold (which is close to ensuring the continuity of the macrodisplacements, see the A-bars in Figure \ref{fig:ddr_flam}). The continuity of the macrodisplacements ($\kmEoM \to \infty$) seemed to be a necessary condition for convergence in the case of large displacements (see the B-bars in Figure \ref{fig:ddr_flam}). With this macroscopic continuity, the anisotropic differentiation of the search directions improved the convergence rate (see the C-bars in Figure \ref{fig:ddr_flam}).}

\vspace{5mm}
{ \noindent \textit{Remark:}
The satisfaction of the macroscopic problem in the admissibility stage (Eq. \eqref{eq:pb_macro}) plays a major role in the handling of global geometric nonlinearities thanks to the transmission of the large-wavelength part of the solution. This requires the expensive continuous updating and assembling of the macroscopic homogenized operator ${}^i\mathbf{L^M}$, which depends on the current configuration of the substructures. There have been some unsuccessful attempts to update the macroscopic problem only at the beginning of each time step, or even only after each local stage, but they resulted in divergence problems or erroneous solutions depending on the time step discretizations. In our work, we used continuous updating, but were able to limit the updating and, thus, the computation time thanks to the introduction of some additional criteria.
}

\subsection{Contact with possible opening}
\label{section:ddr_contact}

Usually, in contact problems, no micro/macro separation of the search directions of the admissibility stage is applied ($\kmEoM = \kmEom = \kmEo$) and the value $\kpEo = \kmEo = E/L_\GammaEo$ seems to be a good choice  \cite{Ladeveze02a}. In the case of open interfaces in slender structures, the best values are $\kpEo = \kmEo =0$. The mechanical behavior of the structure is so different with and without contact that an incorrect setting of the search directions can affect the convergence of the method dramatically. Therefore, we had to develop a search direction updating strategy which takes into account the predicted status of the interface. This updating strategy must be a compromise between performance and stability.

Then we applied this updating strategy to an opening contact problem and a closing contact problem, both initialized with an incorrectly predicted interface status, in order to verify whether the updating algorithm leads to the correct solution.

\subsubsection{An opening contact interface:}

\begin{figure}[!b]
       \centering
       \includegraphics[width=.65 \linewidth]{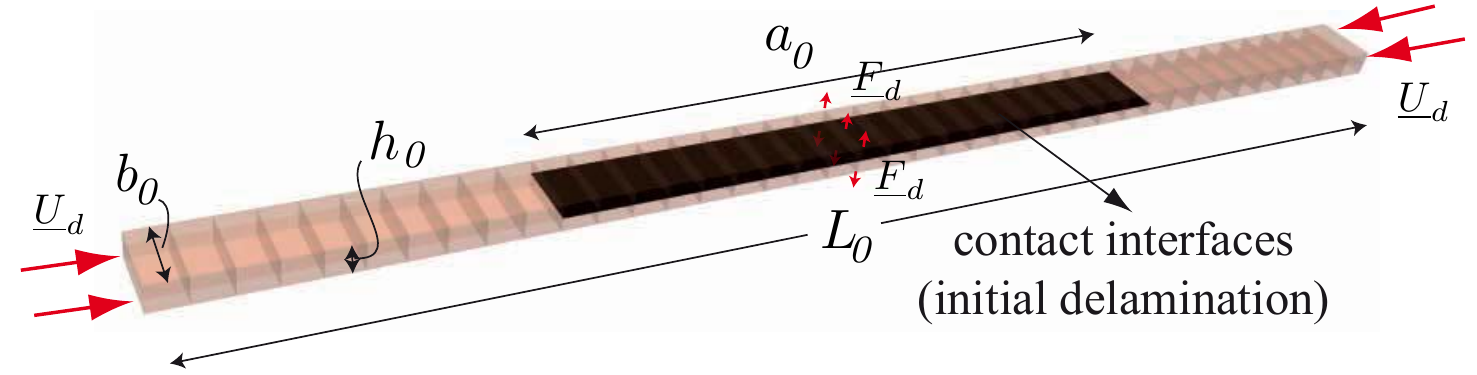}
       \caption{The separation of a contact interface}
       \label{fig:contact_open_confi}
\end{figure}

{The first example concerns a plate built-in at both ends with an initial central delamination $a_\textit0$, subjected progressively to a negative end displacement $\underline U_d$ and to a perturbation consisting of a symmetrical central force $\underline F_d$, as illustrated in Figure \ref{fig:contact_open_confi}}. The data are: $L_\textit0 = 20$ mm, $h_\textit0 = 0.2$ mm, $b_\textit0 = 1$ mm and $a_\textit0 = 10$ mm. The geometry was divided into 80 subdomains and 126 interfaces (see Figure \ref{fig:contact_open_confi}).

\begin{figure}[t]
       \centering
       \includegraphics[width=.65 \linewidth]{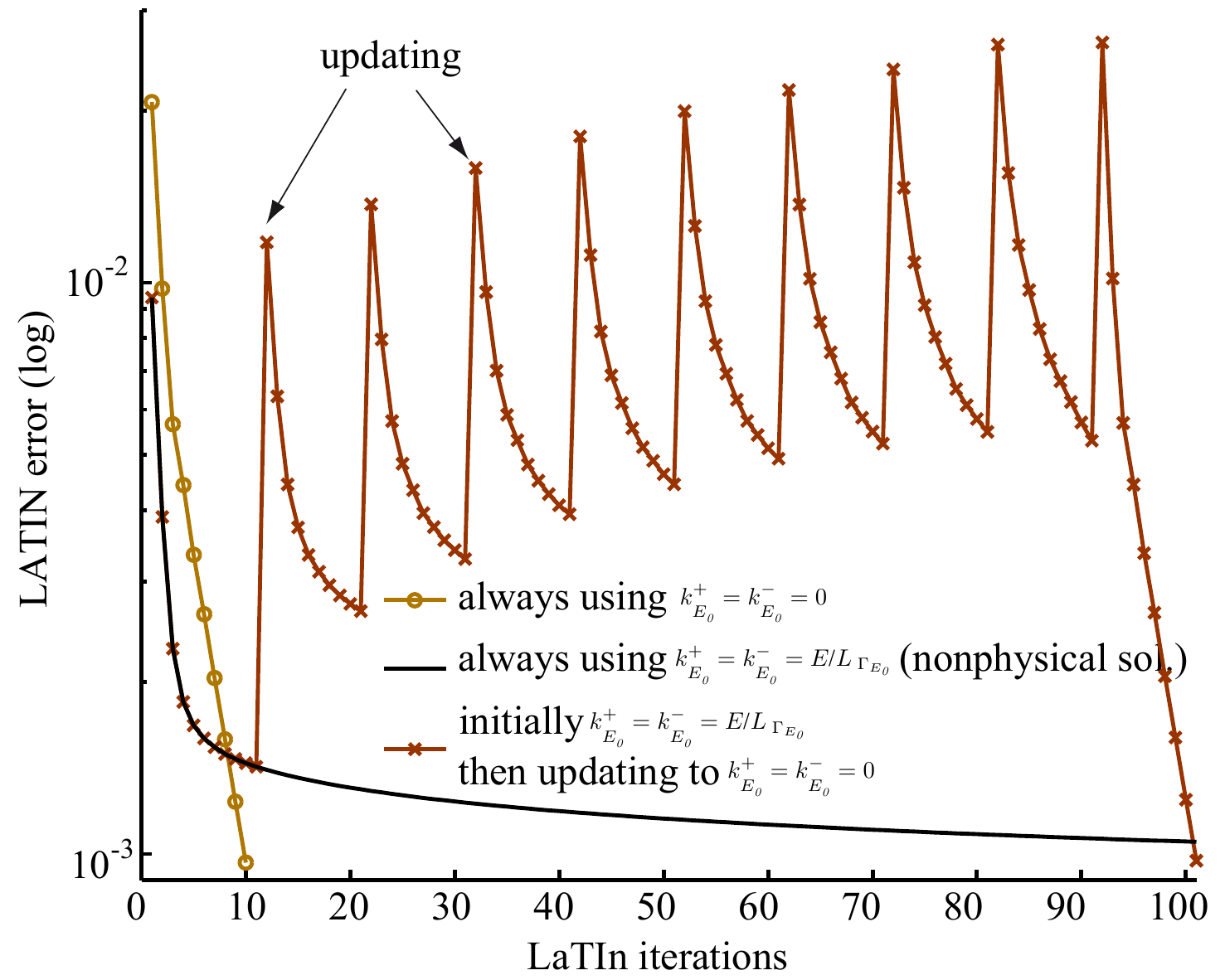}
       \caption{The LATIN error in the open contact (first time step) for different search directions}
       \label{fig:contact_open_ddr}
\end{figure}

\begin{figure}[t]
       \centering
       \includegraphics[width=1 \linewidth]{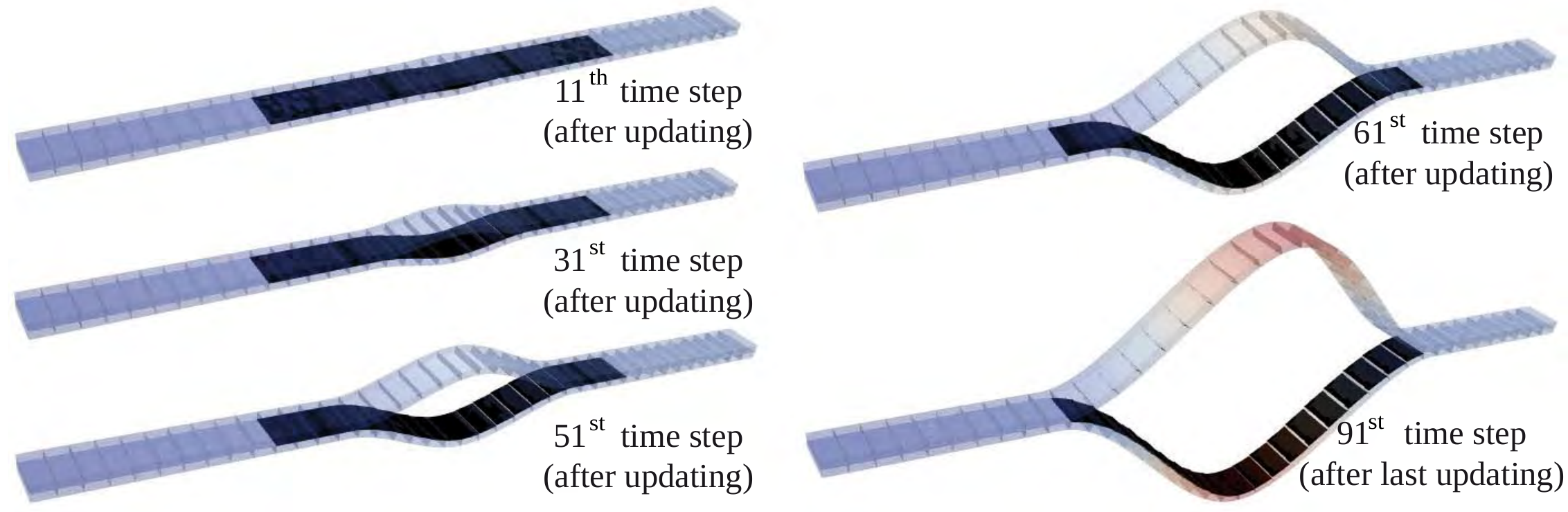}
       \caption{Deformed configurations after updating (magnification factor x500)}
       \label{fig:contact_open_iter}
\end{figure}

Figure~\ref{fig:contact_open_ddr} shows the evolution of the error in the first time step as a function of the number of LATIN iterations. Using the correct guess $\kpEo = \kmEo  = 0$, convergence was achieved in 10 iterations. Conversely, using incorrect values $\kpEo = \kmEo = E/L_\GammaEo$, it was extremely difficult to obtain the correct solution: because of the additional interface stiffness, stagnation in a non-physical configuration occurred. Therefore, we undertook to check the physical status of the interface every 10 iterations and to update the search directions accordingly. Figure~\ref{fig:contact_open_iter} shows different states obtained thanks to this updating strategy.

\subsubsection{A closing contact interface:}

\begin{figure}[!b]
       \centering
       \includegraphics[width=.65 \linewidth]{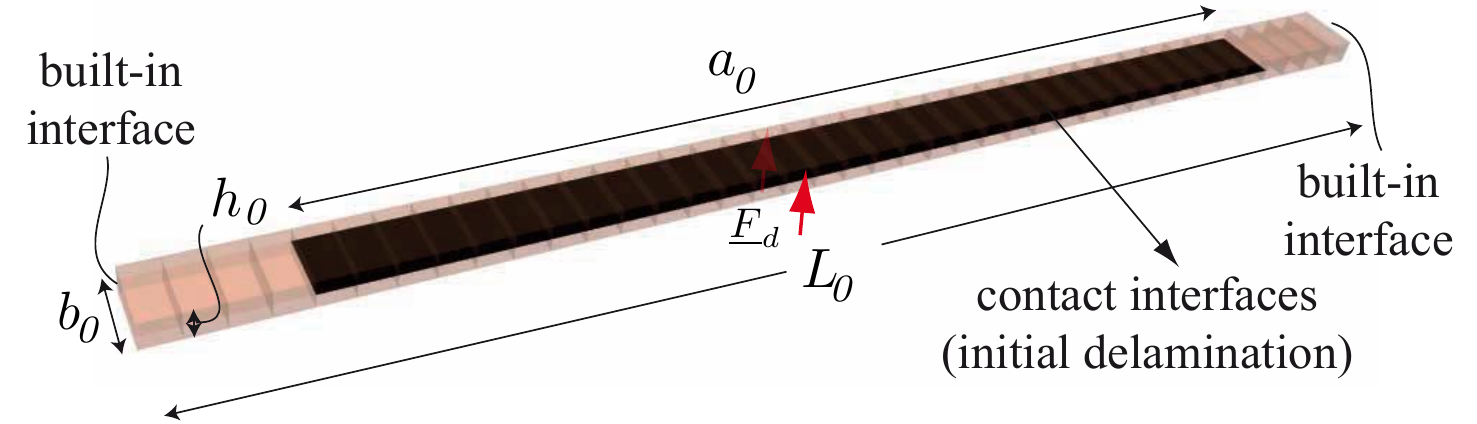}
       \caption{The closing contact interface}
       \label{fig:contact_close_confi}
\end{figure}

{The second example concerns a plate built-in at both ends with an initial central delamination $a_\textit0$ subjected to a central vertical bending force $\underline F_d$ applied to the lower ply, as illustrated in Figure \ref{fig:contact_close_confi}.} The data are: $L_\textit0 = 20$ mm, $h_\textit0 = 0.2$ mm, $b_\textit0 = 1$ mm and $a_\textit0 = 16$ mm. The geometry was divided into 80 substructures and 126 interfaces (see Figure \ref{fig:contact_close_confi}).

\begin{figure}[t]
       \centering
       \includegraphics[width=.65 \linewidth]{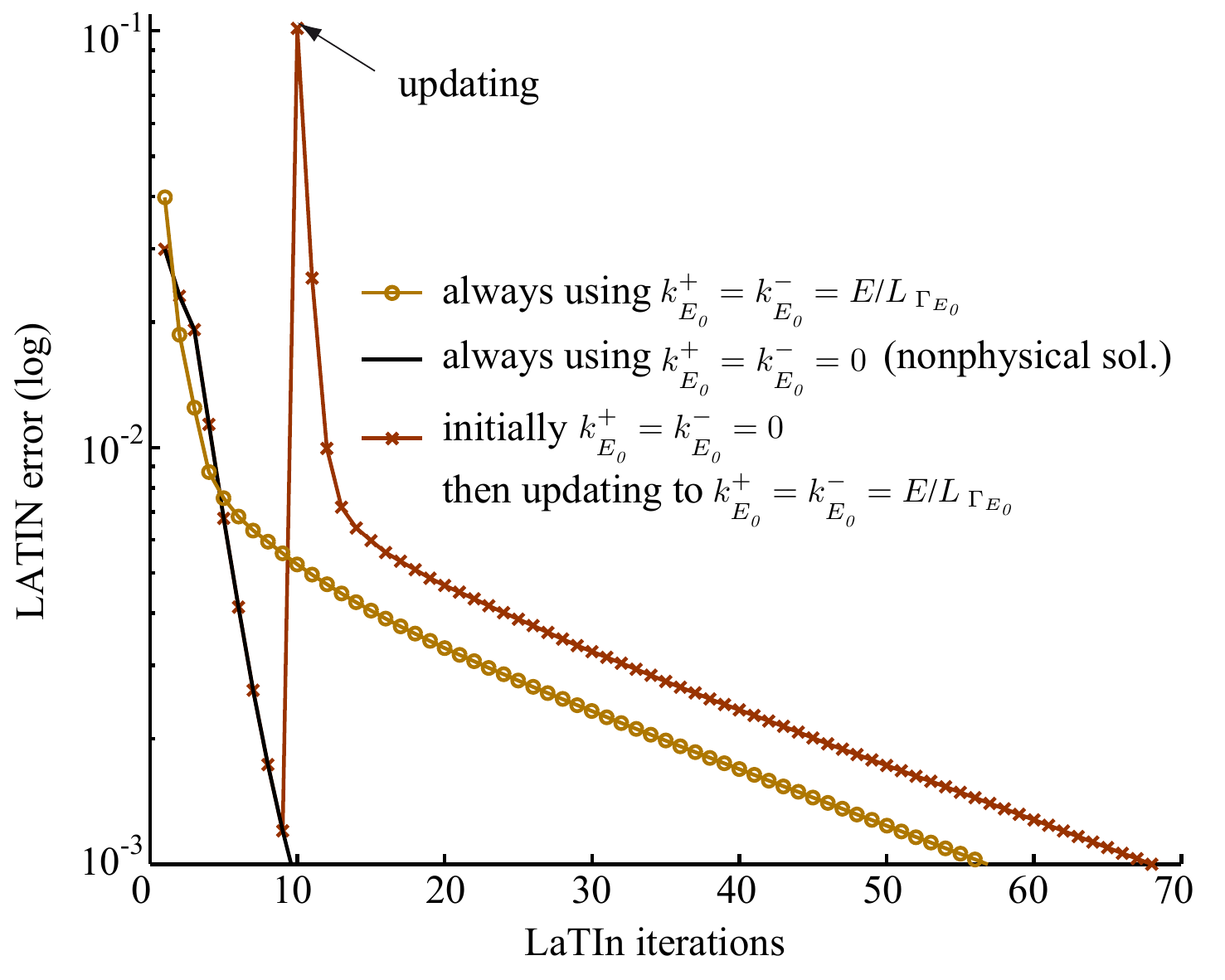}
       \caption{The LATIN error in the closing contact (first time step) for different search directions}
       \label{fig:contact_close_ddr}
\end{figure}

\begin{figure}[t]
       \centering
       \includegraphics[width=.65 \linewidth]{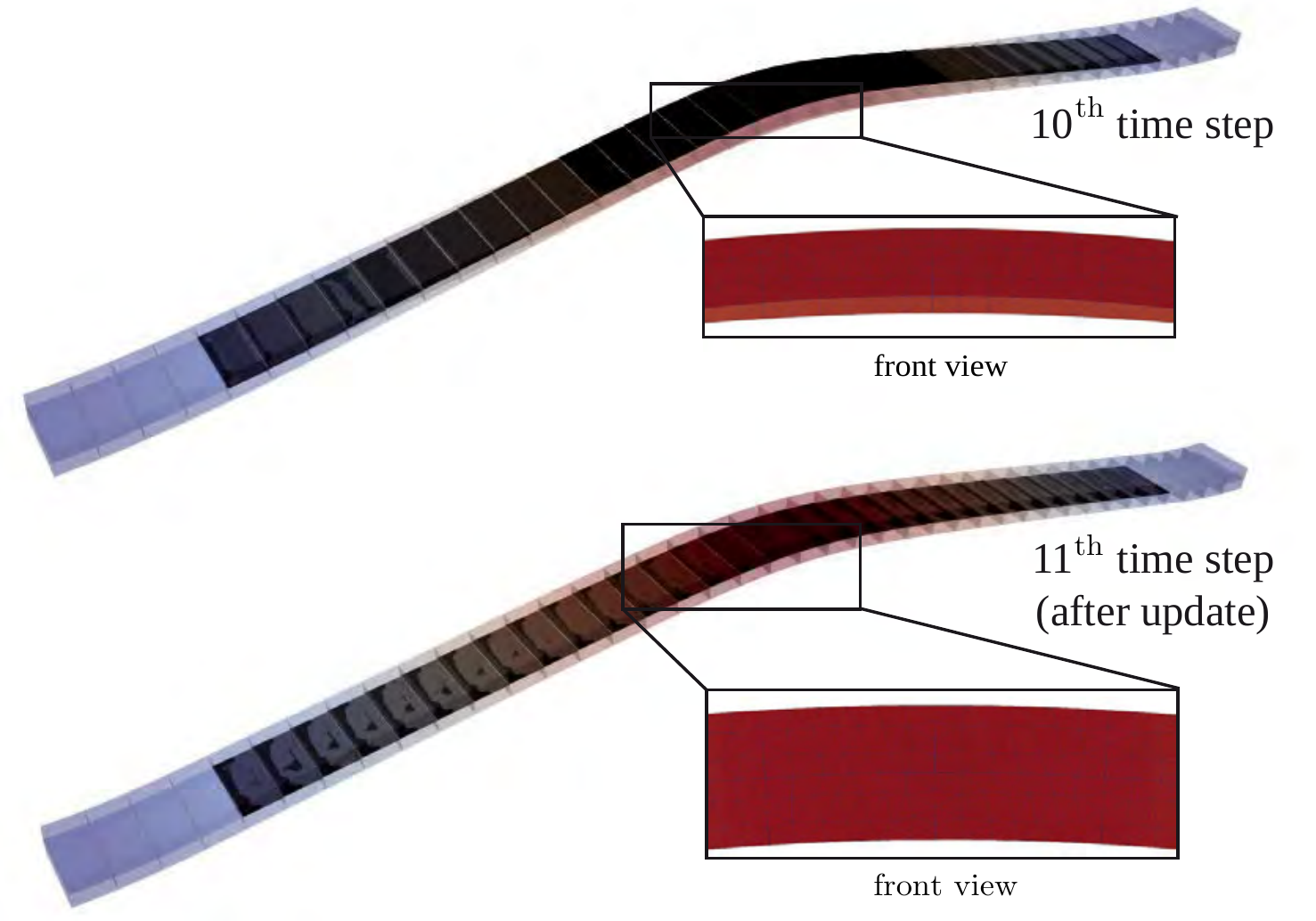}
       \caption{Deformed configuration after updating}
       \label{fig:contact_close_iter}
\end{figure}

In that configuration, the optimum guess for initializing search direction $\kpEo = \kmEo = E/L_\GammaEo$ led to a reasonable number of iterations (see Figure~\ref{fig:contact_close_ddr}). The introduction of an incorrect stiffness (assuming an open interface $\kpEo = \kmEo = 0$) resulted in an erroneous configuration with penetration. It was necessary to update the status in order to obtain a suitable value of the search direction and evaluate the error criterion correctly. Figure~\ref{fig:contact_close_iter} shows the deformed configuration before and after updating. One can observe that overlapping of the plies occurred at the beginning, but was eliminated after updating.

\textit{Remark:} The use of search directions $\kpEo = \kmEo = E/L_\GammaEo$ for the closed interfaces and $\kpEo = \kmEo  = 0$ for the separated interfaces resulted in a proper macroscopic problem (representing the stiffness of the contact interface) and a correct converged solution. The proposed updating scheme according to the interface's state appeared to work properly for an interface going from a closed state to an open state or vice versa. However, in the case of more complex problems, the strategy was not always found to converge because it is difficult to find the exact setup of the updating algorithm in order to avoid divergence. After some empirical tests, we decided to use for both states a unified search direction equal to $\kpEo = \kmEo = (E / L_\GammaEo )/(L_\textit0/h_\textit0)^2$, where $L_\textit0 / h_\textit0$ is the slenderness coefficient of a ply. This choice leads to physically sound solutions with only a few more iterations than using the optimum value for each state of the interface.

\subsection{Cohesive interfaces}
\label{section:ddr_cohesive_interfaces}

Cohesive interfaces lead to the same difficulties as contact interfaces: the recommended values are $\kpEo \to \infty$ and  $\kmEoM = \kmEom = \kmEo = 2 k^0 (1-d)$, and the main problem resides in the definition of a stable and efficient updating strategy for $\kmEo$ with respect to the evolution of $d$.

\begin{figure}[ht]
       \centering
       \includegraphics[width=.65 \linewidth]{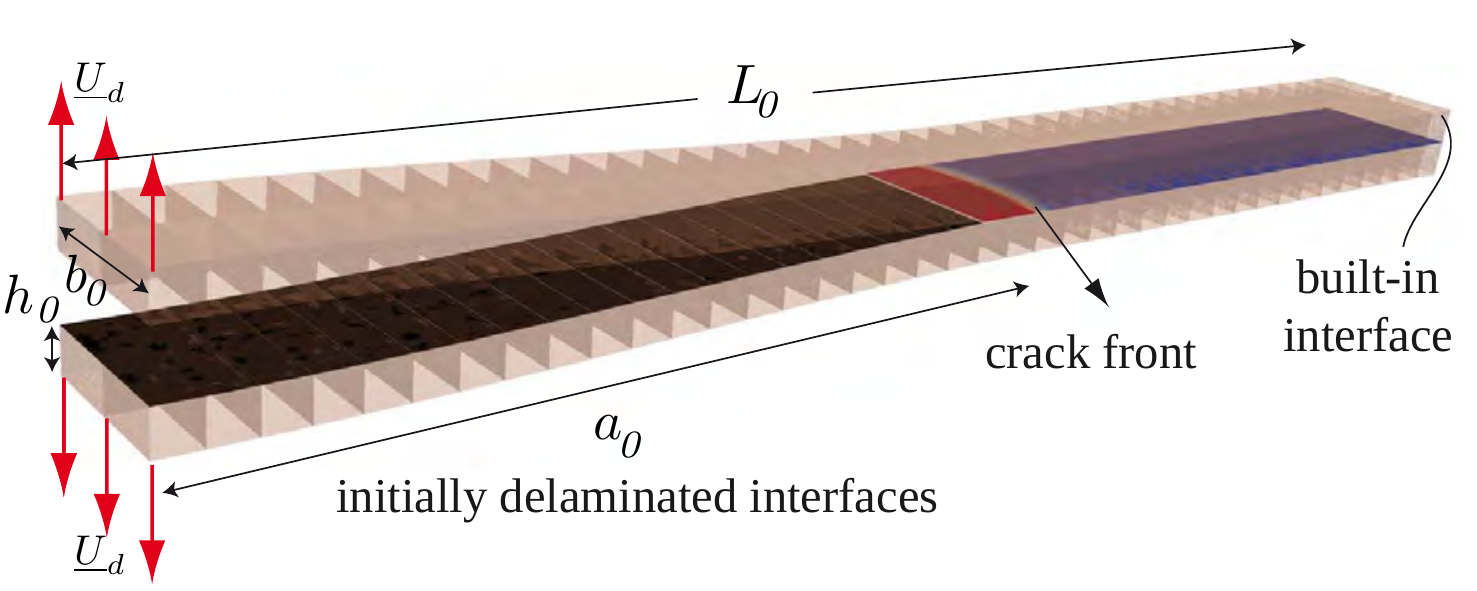}
       \caption{The substructures and interfaces of the deformed DCB specimen}
       \label{fig:delam_confi}
\end{figure}

The strategies were evaluated based on a DCB test under small perturbations. {Figure \ref{fig:delam_confi} shows the deformed specimen with an original crack $a_\textit0$ subjected to a vertical displacement $\underline U_d$ separating the two arms formed by the crack and built-in at the opposite end}. The data are: $L_\textit0=20$~mm, $h_\textit0=0.5$~mm, $b_\textit0=2$~mm, $a_\textit0=10$~mm, $E=135,000$~MPa, $\nu=0.3$ and $k^\textit0=100,000$~N/mm$^3$, $\alpha=1$, $n=0.5$ and $Y_c=0.4$~N/mm. The structure was divided into 160 substructures and 324 interfaces (see Figure \ref{fig:delam_confi}), which amounted to 550,000 DOFs for the whole mesh and 2,916 DOFs for the macroscopic problem. Because of its small size, the macroscopic problem was solved using a direct solver. There were at least 10 elements in the process zone and the load was applied in $16$ time steps.

\begin{figure}[ht]
       \centering
       \includegraphics[width=.65 \linewidth]{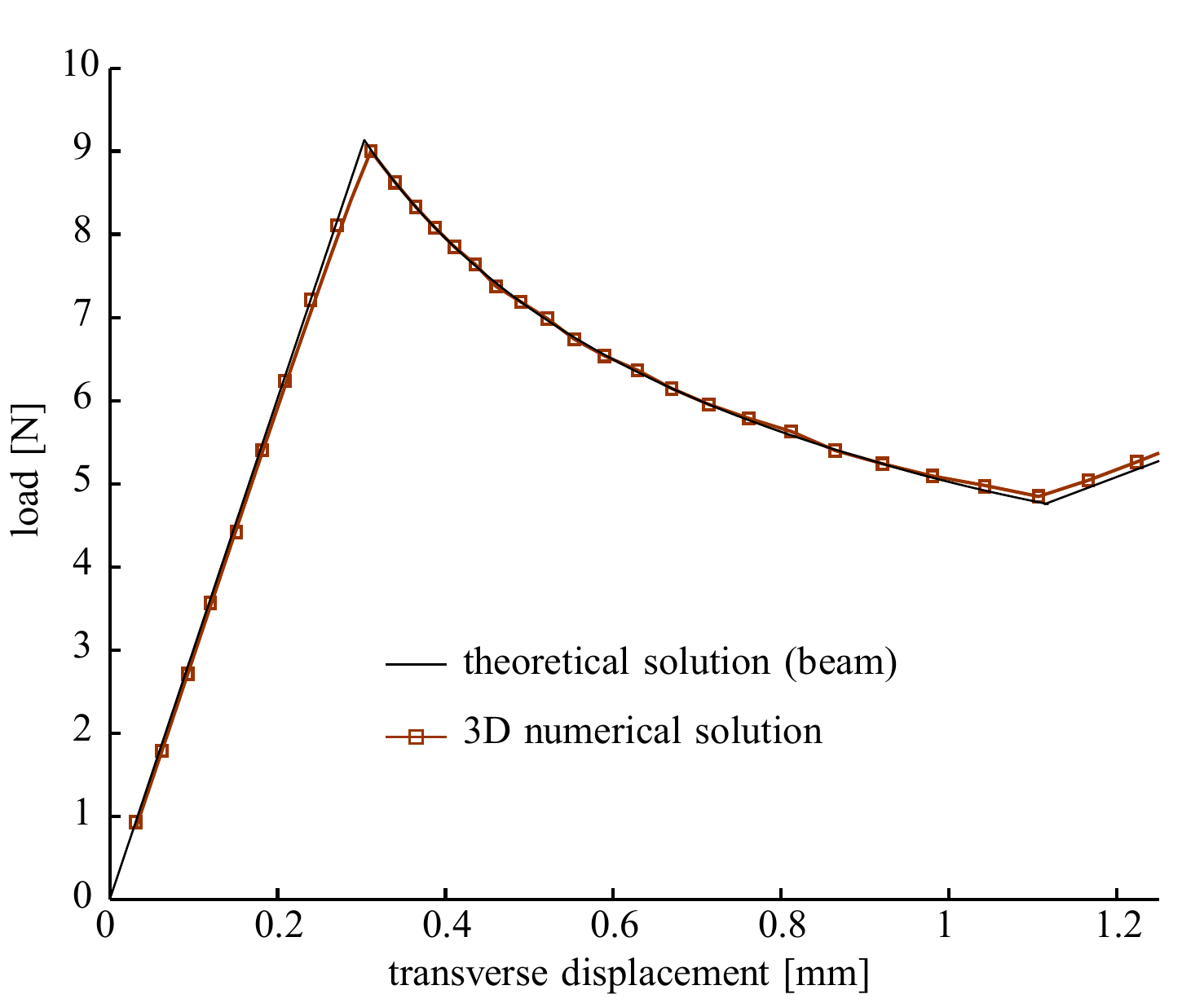}
       \caption{The load-displacement curve of the DCB test}
       \label{fig:delam}
\end{figure}

{Figure \ref{fig:delam} shows the load-displacement curve at the end of one of the two arms of the specimen. The global response until the complete failure of the specimen agrees with the theoretical solution calculated for a beam on an elastic foundation \cite{Kanninen73} using linear fracture mechanics theory for the propagation \cite{Allix95}. The difference in slope which can be observed before the softening part of the curve is due to the fact that the damage law used in the calculation enables damage to take place progressively before the rupture of a point. For the comparison of the different updating procedures, only the first $16$ time steps were considered.}

\begin{figure}[ht]
       \centering
       \includegraphics[width=.65 \linewidth]{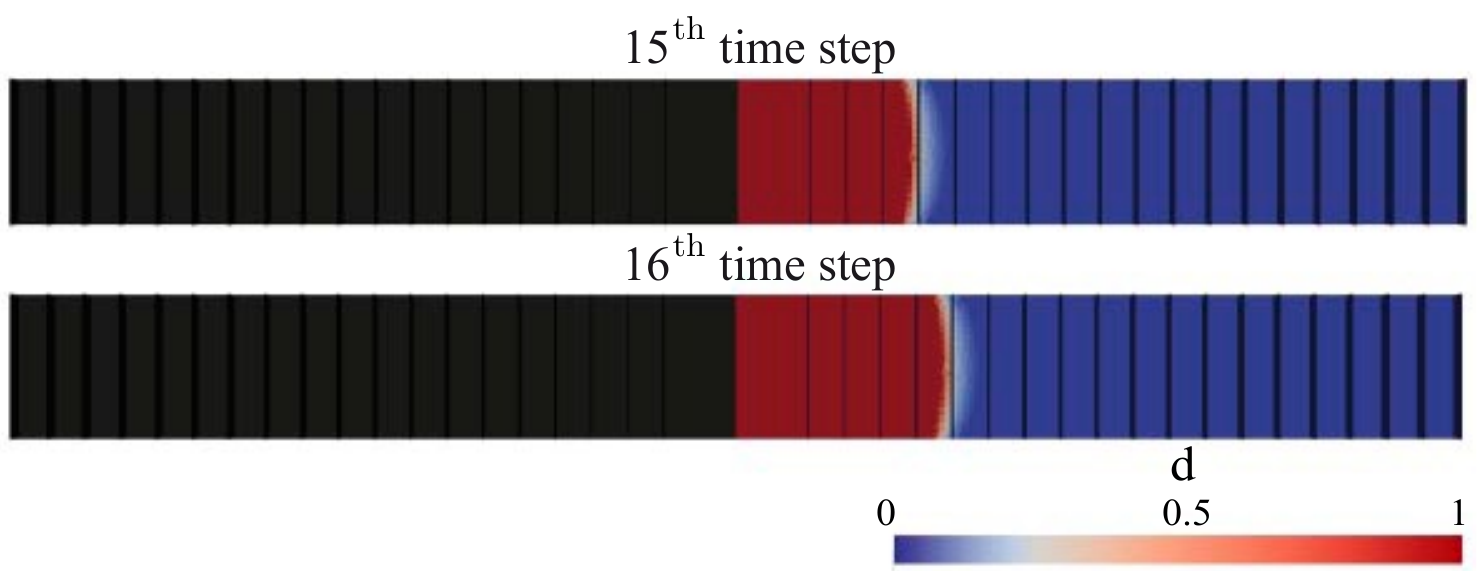}
       \caption{The crack front and the interfaces after the $15^{th}$ and $16^{th}$ time steps}
       \label{fig:crack_front}
\end{figure}

{Figure \ref{fig:crack_front} shows the crack front after the $15^{th}$ and $16^{th}$ time steps.} It is important to note that a first element becomes fully damaged after the $10^{th}$ time step; subsequently, the process zone expands to a new subdomain with each new time step to reach five completely damaged interfaces after the $16^{th}$ step.

\begin{figure}[t]
       \centering
       \includegraphics[width=.65 \linewidth]{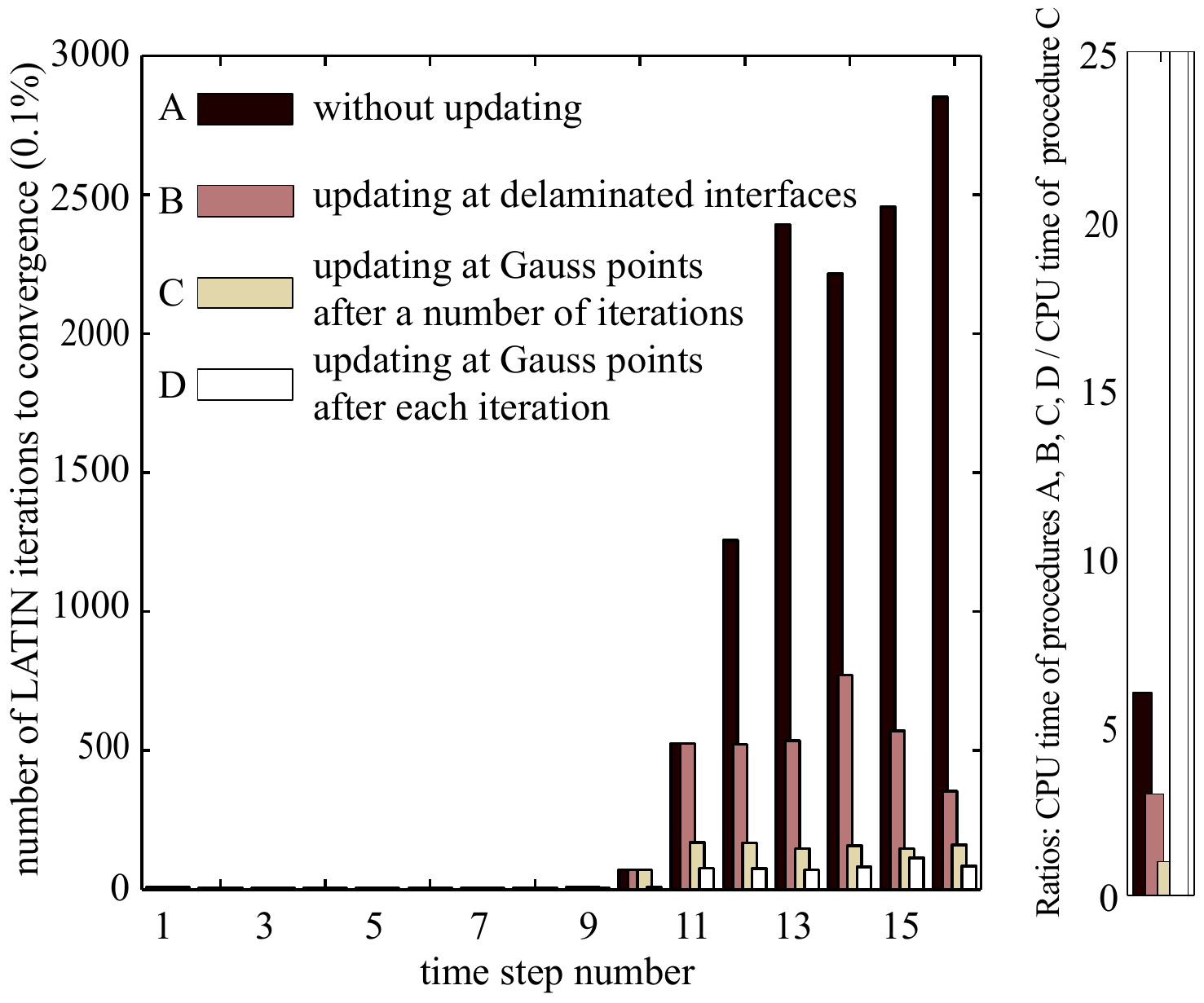}
       \caption{The convergence rate and the corresponding CPU time}
       \label{fig:ddr_delam}
\end{figure}

Figure \ref{fig:ddr_delam} shows the convergence rate of each time step and the corresponding total CPU time for the different updating schemes. As expected, the number of iterations increases after the $10^{th}$ step.Without updating, the number of iterations explodes (A-bars), although the CPU time remains moderate. If, as proposed in \cite{Kerfriden09}, updating is performed only at the interfaces which have been fully damaged (B-bars), both the number of iterations and the CPU time decrease respectively by a factor of 3 and a factor of 2. The optimum strategy seems to consist in updating at each Gauss point after a predetermined number of LATIN iterations (\emph{e.g.}  100), in which case both the number of iterations and the CPU time decrease another $65\%$ compared to the previous strategy and become very small. The more aggressive strategy which consists in updating at each Gauss point after each iteration (D-bars) leads to the smallest number of iterations, but a huge CPU time (25 times that of the procedure C).

\section{Examples of combined buckling and delamination}
\label{section:flam_delam}

\subsection{A built-in beam under compressive loading}

\begin{figure}[!t]
       \centering
       \includegraphics[width=.65 \linewidth]{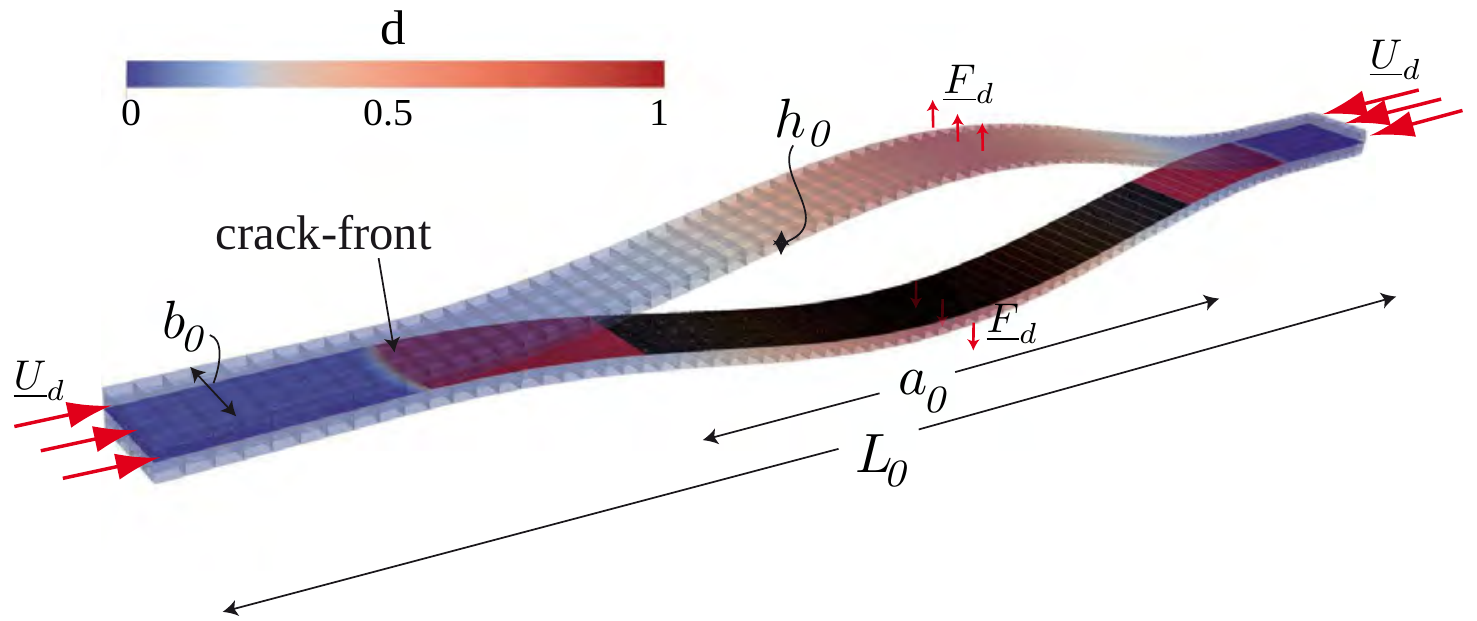}
       \caption{A compressed plate with an initial delamination in local buckling: substructures, deformed configuration (with no amplification) and crack front}
       \label{fig:flam_delam_confi}
\end{figure}

The example discussed here consists in the simulation of a built-in plate with a central initial delamination $a_\textit0$ subjected to an axial compressive loading $\underline{U}_d$ and a symmetrical central perturbation $\underline{F}_d$ (see Figure \ref{fig:flam_delam_confi}). The data are:
$L_\textit0 = 20$ mm, $h_\textit0 = 0.2$ mm, $b_\textit0 = 1$ mm, $a_\textit0=10$~mm, $E = 135,000$ MPa, $\nu = 0.3$, $k^\textit0_n = k^\textit0_t = 100,000$ N/mm$^3$, $\alpha = 1$, $n=0.5$ and $Y_c = 0.4$ N/mm. The geometry was divided into 1,280 substructures and 3,248 interfaces (see Figure \ref{fig:flam_delam_confi}), leading to a mesh totaling 1.4 million DOFs with 24 linear wedge elements through the thickness (12 elements in each ply). The macroscopic problem and the supermacroscopic problem contained respectively 29,232 DOFs and 132 DOFs, and were solved using 24 processors.

The solution of this problem is characterized by the competition between global and local buckling. The critical local buckling load is approximately twice that of a beam of length $a_\textit0$ built-in at both ends:
\begin{equation*}
P_{c}^{local} = 2 \left(\frac{4 \pi ^2 E I_\textit0}{a_\textit0 ^2}\right) =71.06 \; \text{N}\;,
\end{equation*}
where $I_\textit0$ is the second moment of area of a single ply. The global buckling load of a beam with an initial crack is given approximately by the formula \cite{Allix99}:
\begin{equation*}
P_{c}^{global} =  \left( \frac{L_\textit0-a_\textit0}{L_\textit0}\right) \left(\frac{4 \pi ^2 E 8I_\textit0}{L_\textit0 ^2}\right) + \left( \frac{2a_\textit0}{L_\textit0}\right) \left(\frac{4 \pi ^2 E I_\textit0}{L_\textit0 ^2}\right)= 44.41\; \text{N}\ .
\end{equation*}

{In the case of local buckling without imperfection, the propagation condition can be approximated by the formula \cite{Bruno90}:
\begin{equation*}
\frac{3}{16}\xi^4 + 2\xi^2 = \frac{4G_cb_ \textit0}{\pi^2P_{c}^{local}}\;,
\end{equation*}
where $\xi = w_{L_ \textit0 /2}/ a_ \textit0$ is the dimensionless transverse displacement parameter and $G_c$ is the critical energy release rate (for the damage interface law used here, $G_c = Y_c$). Thus, the propagation condition in local buckling is:
\begin{equation}
w_{L_ \textit0 /2} = 0.34 \; \text{mm}\;. \label{delam_condition}
\end{equation}
}

\begin{figure}[!h]
       \centering
       \includegraphics[width=.65 \linewidth]{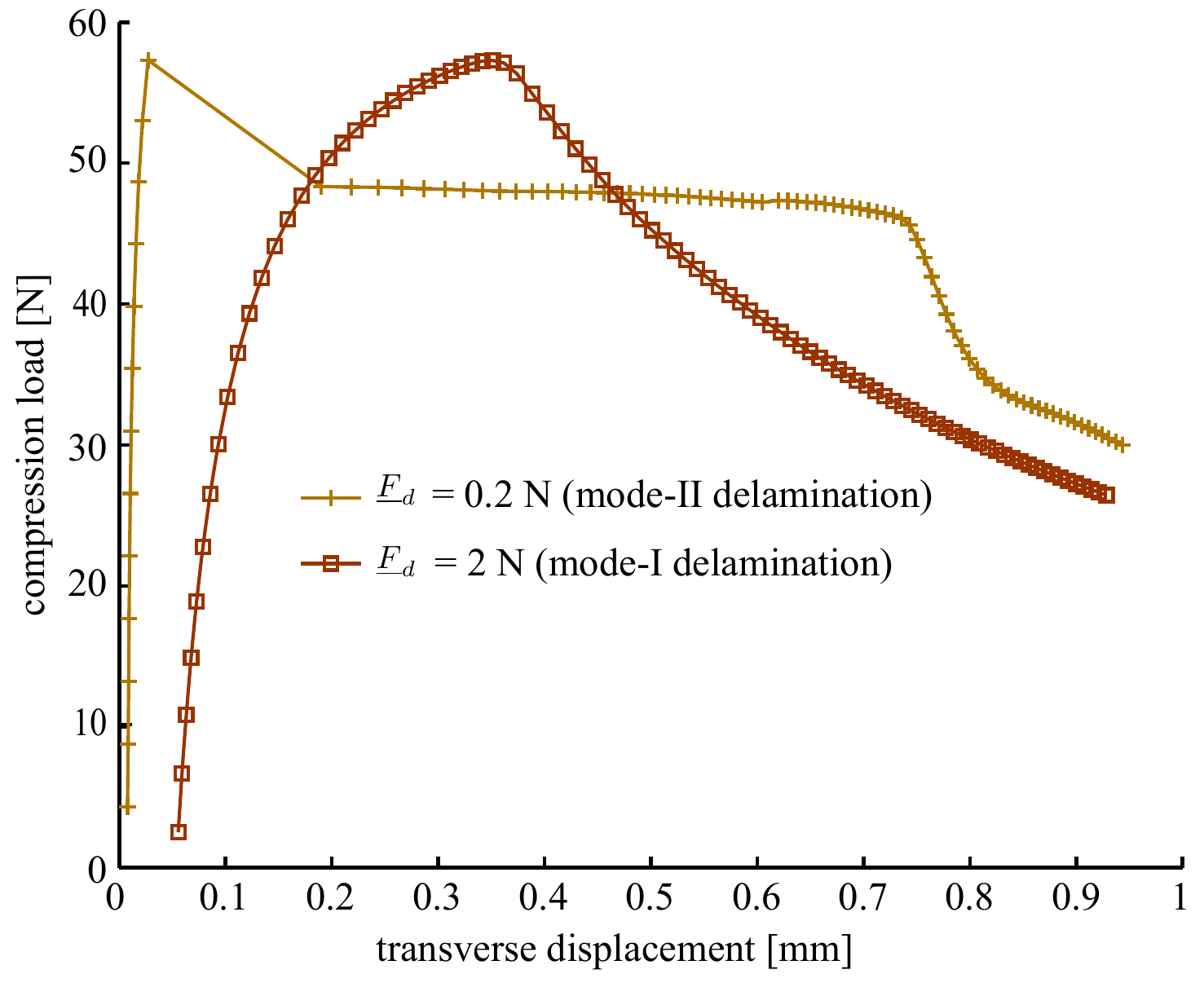}
       \caption{The load-displacement curve of the compressed plate with an initial delamination}
       \label{fig:flam_delam}
\end{figure}

\begin{figure}[t]
       \centering
       \includegraphics[width=.65 \linewidth]{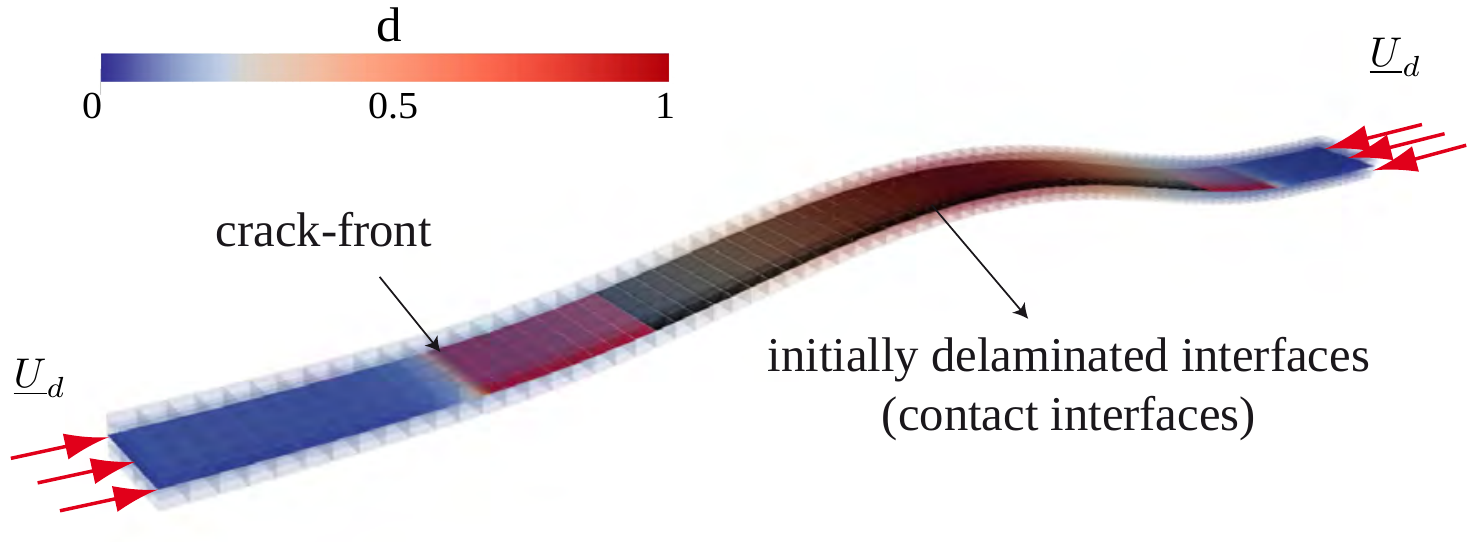}
       \caption{The compressed plate with an initial delamination in global buckling: substructures, deformed solution in global buckling (without amplification) and crack front}
       \label{fig:flam_delam_confi2}
\end{figure}

First, we studied the numerical response of the structure under a small symmetrical perturbation ($\underline{F}_d=0.2~N$). From the evolution of the compressive load as a function of the maximum transverse displacement of the upper ply shown in Figure \ref{fig:flam_delam}, one can see that before the load reaches the critical local value $P_{c}^{local}=71.06$~N the response switches from local buckling to global buckling (after the $14^{th}$ time step), which reduces the load to about the critical global value $P_{c}^{global} = 44.41$~N and induces mode-II delamination at about $0.75$~mm maximum transverse displacement. The deformed configuration and the corresponding crack front after the last time step are shown in Figure \ref{fig:flam_delam_confi2}. An increase in the amplitude of the symmetrical perturbation to $\underline{F}_d=2$~N is enough to determine the local buckling mode. Figure \ref{fig:flam_delam} indicates a local buckling load equal to about $56$~N, which is less than the approximate evaluation. Mode-I delamination begins at about $0.36$~mm maximum transverse displacement, {which is very close to the theoretical delamination condition \eqref{delam_condition}}. The deformed configuration and its crack front after the last time step are shown in Figure \ref{fig:flam_delam_confi}. {The solution obtained can be compared to the numerical solution of a beam given in \cite{Allix99}. There is good agreement in terms of critical load and delamination propagation.}

\subsection{Multiple through-the-width delaminations}
\label{section:multi_delaminations}

\begin{figure}[!b]
       \centering
       \includegraphics[width=.65 \linewidth]{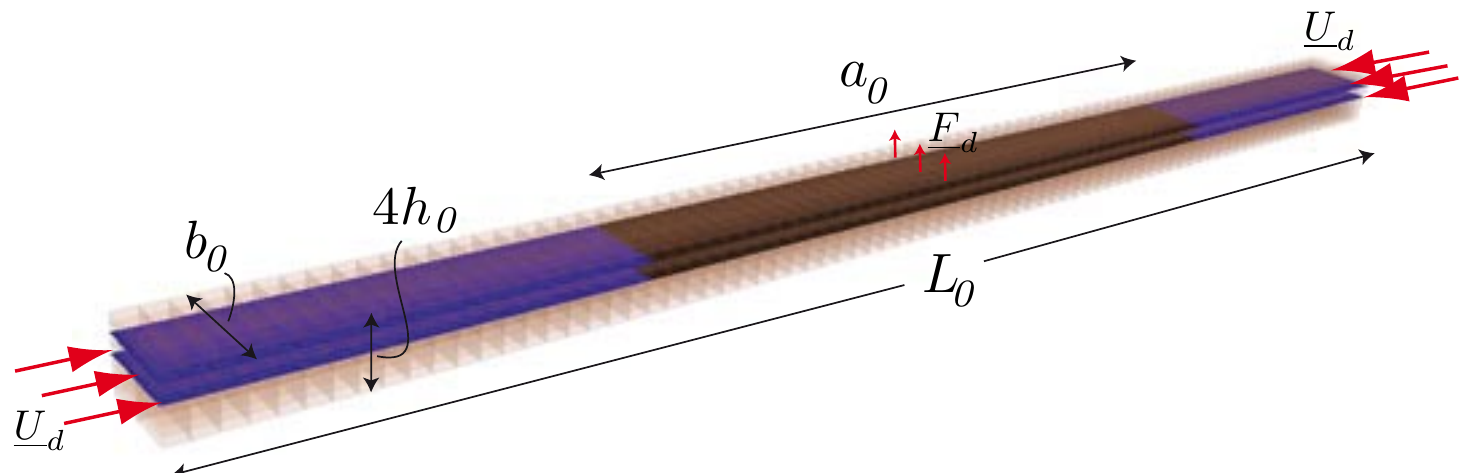}
       \caption{A compressed laminate with two initial delaminations}
       \label{fig:flam_delam_asym2_confi}
\end{figure}

{This section deals with a built-in laminate made of 4 plies of thickness $h_\textit0$ with two initial central delaminations $a_\textit0$ at $2h_\textit0$ and $3h_\textit0$, subjected progressively to a negative end displacement $\underline U_d$ and a central perturbation $\underline F_d$ applied to the upper ply as shown in Figure \ref{fig:flam_delam_asym2_confi}.} The data are: $L_\textit0 = 20$~mm, $h_\textit0=0.2$~mm, $b_\textit0=1$~mm, $a_\textit0=10$~mm, $E_1=185,500$~MPa, $E_2=E_3=9,900$~MPa, $\nu_{12}=\nu_{13}=0.34$, $\nu_{23}=0.5$, $G_{12}=G_{13}=6,160$~MPa, $G_{23}=3,080$~MPa, $k^\textit0_n=k^\textit0_t=100,000$~N/mm$^3$, $\alpha=1$, $n=0.5$ and $Y_c=0.4$~N/mm. The lay-up sequence is $[0¡/90¡]_s$. The geometry was divided into 1,280 substructures and 3,064 interfaces (see Figure \ref{fig:flam_delam_asym2_confi}), leading to a mesh totaling 2 million DOFs with 12 linear wedge elements through the thickness of each ply. The macroscopic problem and the supermacroscopic problem contained respectively 27,576 DOFs and 168 DOFs, and were solved using 30 processors.

\begin{figure}[!t]
       \centering
       \includegraphics[width=.65 \linewidth]{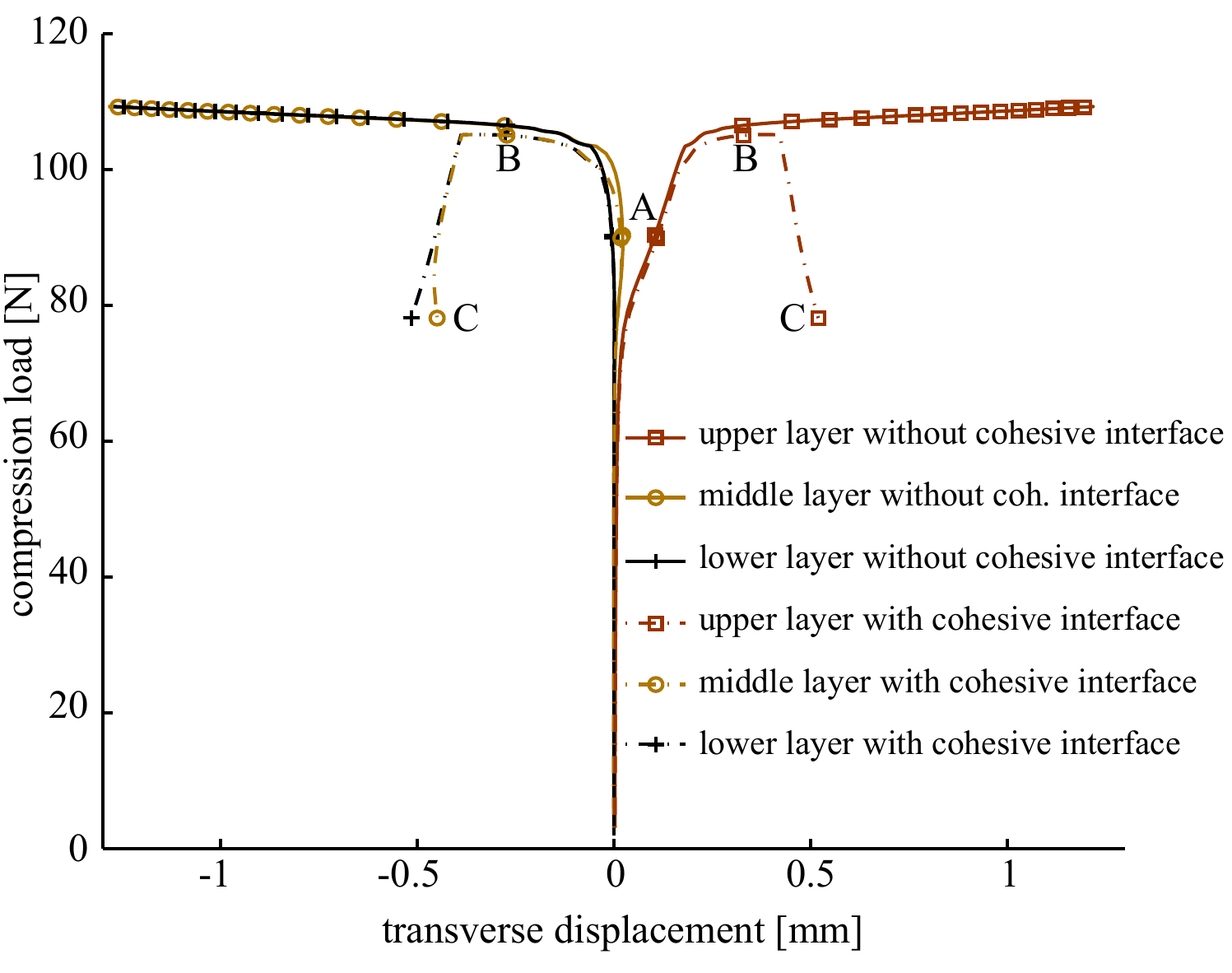}
       \caption{The load-displacement curves of the three layers of a compressed laminate with two initial delaminations. Note: only some time steps are represented by markers.}
       \label{fig:flam_delam_asym2}
\end{figure}

A first computation without cohesive interfaces was carried out in order to calculate the post-buckling response of the laminate. In Figure \ref{fig:flam_delam_asym2}, the solid lines represent the compressive load as a function of the transverse displacement in the middle of each of the three layers determined in the laminate by the two initial cracks. One can see that, first, the upper layer undergoes nonsymmetrical local buckling after $80$ N. Then, at $100$ N, the middle and lower layers (\emph{i.e.} the unperturbed layers) buckle in the opposite direction, leading to a symmetrical local buckled configuration. For the nonsymmetrical buckled shape, the transverse displacement of the middle layer is positive (pulled by the upper ply) and is separated from the lower layer (point A). Subsequently, in the symmetrically buckled shape, the transverse displacement of the middle ply becomes negative (pulled by the buckling of the lower layer) and the lower and middle layers are in contact (point B). Another calculation without initial delamination resulted in a global buckling load more than five times the critical load with cracks.

\begin{figure}[ht]
       \centering
       \includegraphics[width=.65 \linewidth]{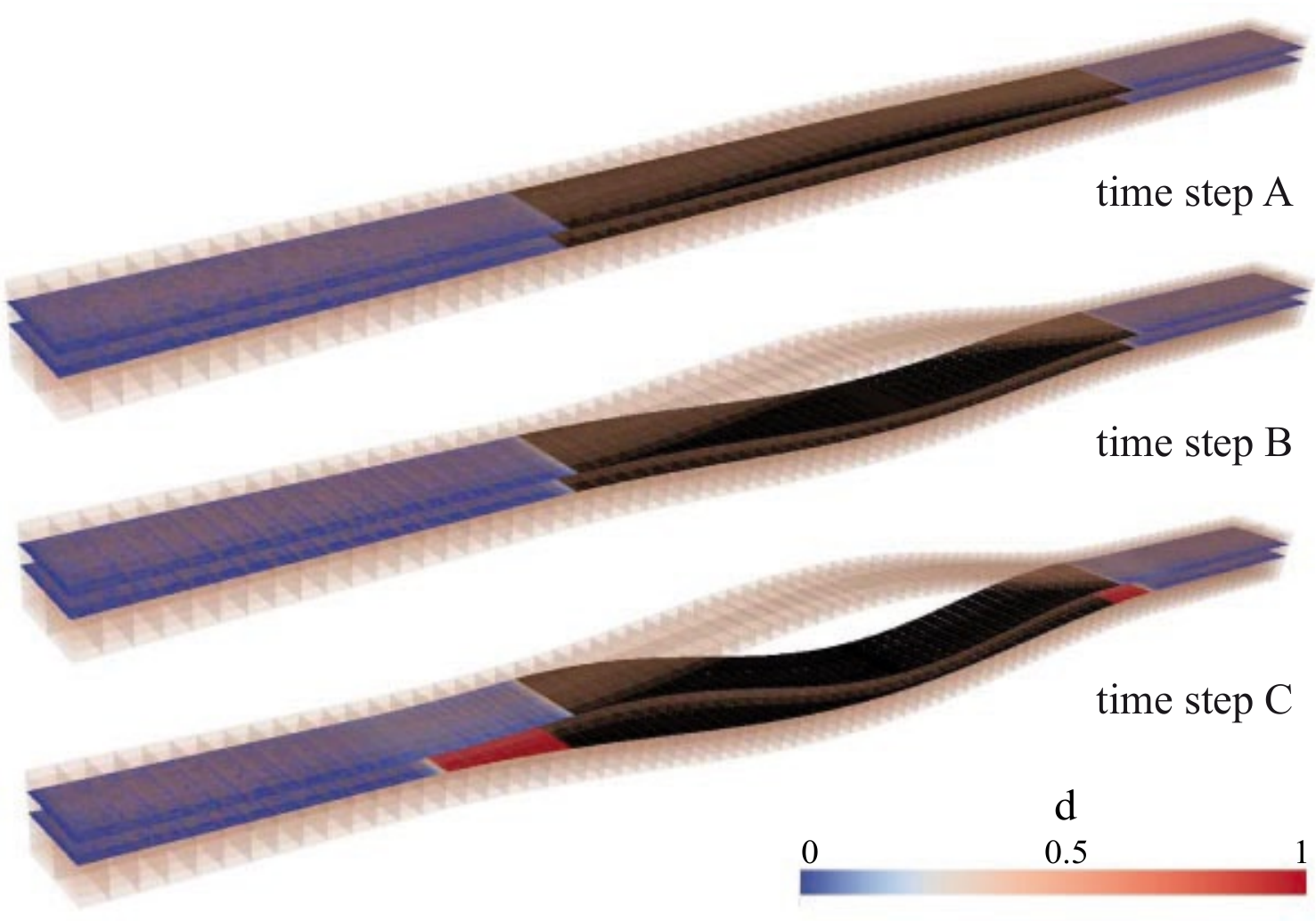}
       \caption{The deformed configurations and crack fronts after three time steps of a compressed laminate with two initial delaminations}
       \label{fig:flam_delam_asym2_iter}
\end{figure}

The response taking into account delamination between plies is also shown in Figure \ref{fig:flam_delam_asym2}. The buckling behavior is almost the same as the response without cohesive interfaces. The critical buckling load is smaller because of the damage interface law. After point B, the load decreases due to the fact that the first elements in the crack front are completely damaged. At point C, the middle and lower layers are no longer in contact. The deformed configurations and crack fronts after time steps A, B and C of Figure \ref{fig:flam_delam_asym2} are shown in Figure \ref{fig:flam_delam_asym2_iter}.

\section{Conclusion}
\label{sec:conclusion}

We presented an efficient calculation of extensive delamination in the presence of geometrically nonlinear effects thanks to a three-scale domain decomposition strategy based on an iterative algorithm. This method enables the resolution of huge nonlinear systems of equations on the most suitable scales.

The LATIN algorithm we proposed can handle both the geometric nonlinearities in the substructures and the surface degradations at the interfaces, thus making the treatment of both nonlinearities easier. The second scale (\emph{i.e.} the macroscopic problem) enables the rapid transmission of the large-wavelength part of the response through the introduction of an updated homogenized behavior of the substructures which takes into account the deformed configuration and the current damage state of the interfaces. The macroscopic problem was solved using a parallel iterative solver.

The classical values of the method's parameters (\emph{i.e.} the search directions) were shown to be inadequate for slender structures. Therefore, we modified them in order to ensure the scalability of the method, its independence with respect to the geometry of the subdomains, an efficient convergence rate and an adequate CPU time for the treatment of combined buckling and delamination. Examples showing the capabilities of the strategy were also presented.

In subsequent developments, the introduction of large sliding contact conditions in the delaminated area and the coupling of this 3D model with a plate model in the low-gradient zones should be envisaged.

\textbf{Acknowledgement: }  The research leading to these results has received funding from the European Community's Seventh Framework Program FP7/2007-2013 under grant agreement n$^\circ$213371.

\bibliographystyle{plain}
\bibliography{biblio_paper}

\begin{thebibliography}{10}

\bibitem{Allix99}
O.~Allix and A.~Corigliano.
\newblock Geometrical and interfacial non-linearities in the analysis of
  delamination in composites.
\newblock {\em International Journal of Solids and Structures},
  36(15):2189--2216, 1999.

\bibitem{Allix10}
O.~Allix, P.~Kerfriden, and P.~Gosselet.
\newblock On the control of the load increments for a proper description of
  multiple delamination in a domain decomposition framework.
\newblock {\em International Journal for Numerical Methods in Engineering},
  83(11):1518--1540, 2010.

\bibitem{Allix95}
O.~Allix, P.~Ladev{\`e}ze, and A.~Corigliano.
\newblock Damage analysis of interlaminar fracture specimens.
\newblock {\em Composite Structures}, 31(1):61--74, 1995.

\bibitem{Allix98}
O.~Allix, D.~L{\'e}v{\^e}que, and L.~Perret.
\newblock Identification and forecast of delamination in composite laminates by
  an interlaminar interface model.
\newblock {\em Composites Science and Technology}, 58(5):671--678, 1998.

\bibitem{Bottega83}
W.~J. Bottega and A.~Maewal.
\newblock Delamination buckling and growth in laminates.
\newblock {\em Journal of Applied Mechanics - Transactions of The ASME},
  50(1):184--189, 1983.

\bibitem{Bruno00}
D.~Bruno and F.~Greco.
\newblock An asymptotic analysis of delamination buckling and growth in layered
  plates.
\newblock {\em International Journal of Solids and Structures},
  37(43):6239--6276, 2000.

\bibitem{Bruno90}
D.~Bruno and A.~Grimaldi.
\newblock Delamination failure of layered composite plates loaded in
  compression.
\newblock {\em International Journal of Solids and Structures}, 26(3):313--330,
  1990.

\bibitem{Chai81}
H.~Chai, C.~D. Babcock, and W.~G. Knauss.
\newblock One dimensional modeling of failure in laminated plates by
  delamination buckling.
\newblock {\em International Journal of Solids and Structures},
  17(11):1069--1083, 1981.

\bibitem{Champaney99}
L.~Champaney, J.~Y. Cognard, and P.~Lad{\`e}veze.
\newblock Modular analysis of assemblages of three-dimensional structures with
  unilateral contact conditions.
\newblock {\em Computers \& Structures}, 73(1-5):249 -- 266, 1999.

\bibitem{Cochelin91}
B.~Cochelin and M.~Potier-Ferry.
\newblock A numerical model for buckling and growth of delaminations in
  composite laminates.
\newblock {\em Computer Methods in Applied Mechanics and Engineering},
  89(1-3):361--380, 1991.

\bibitem{Cresta07}
P.~Cresta, O.~Allix, C.~Rey, and S.~Guinard.
\newblock Nonlinear localization strategies for domain decomposition methods:
  Application to post-buckling analyses.
\newblock {\em Computer Methods in Applied Mechanics and Engineering},
  196:1436--1446, 2007.

\bibitem{Daridon02}
L~Daridon and K~Zidani.
\newblock {The stabilizing effects of fiber bridges on delamination cracks}.
\newblock {\em Composites Science and Technology}, {62}({1}):{83--90}, {2002}.

\bibitem{Evans84}
A.~G. Evans and J.~W. Hutchinson.
\newblock On the mechanics of delamination and spalling in compressed films.
\newblock {\em International Journal of Solids and Structures}, 20(5):455--466,
  1984.

\bibitem{Guidault08}
P.~A. Guidault, O.~Allix, L.~Champaney, and C.~Cornuault.
\newblock A multiscale extended finite element method for crack propagation.
\newblock {\em Computer Methods in Applied Mechanics and Engineering},
  197(5):381--399, 2008.

\bibitem{Kachanov76}
L.~M. Kachanov.
\newblock Separation failure of composite materials.
\newblock {\em Polymer Mechanics}, 12:812--815, 1976.

\bibitem{Kanninen73}
M.~F. Kanninen.
\newblock An augmented double cantilever beam model for studying crack
  propagation and arrest.
\newblock {\em International Journal of Fracture}, 9(1):83--92, 03 1973.

\bibitem{Kardomateas93}
G.~A. Kardomateas.
\newblock The initial postbuckling and growth-behavior of internal
  delaminations in composite plates.
\newblock {\em Journal of Applied Mechanics - Transactions of The ASME},
  60(4):903--910, Dec 1993.

\bibitem{Kerfriden09}
P.~Kerfriden, O.~Allix, and P.~Gosselet.
\newblock A three-scale domain decomposition method for the {3D} analysis of
  debonding in laminates.
\newblock {\em Computational Mechanics}, 44(3):343--362, 08 2009.

\bibitem{Kosel05}
Franc Kosel, Joze Petrisic, Boris Kuselj, Tadej Kosel, Viktor Sajn, and Mihael
  Brojan.
\newblock Local buckling and debonding problem of a bonded two-layer plate.
\newblock {\em Archive of Applied Mechanics}, 74(10):704--726, 09 2005.

\bibitem{Ladeveze99}
P.~Ladev{\`e}ze.
\newblock {\em Nonlinear Computational Structural Mechanics - {N}ew Approaches
  and Non-incremental Methods of Calculation}.
\newblock Springer-Verlag, Berlin, 1999.

\bibitem{Ladeveze01}
P~Ladev{\`e}ze, O~Loiseau, and D~Dureisseix.
\newblock {A micro-macro and parallel computational strategy for highly
  heterogeneous structures}.
\newblock {\em International Journal for Numerical Methods in Engineering},
  {52}({1-2}):{121--138}, {2001}.

\bibitem{Ladeveze02}
P.~Ladev{\`e}ze and G.~Lubineau.
\newblock An enhanced mesomodel for laminates based on micromechanics.
\newblock {\em Composites Science and Technology}, 62(4):533--541, 2002.

\bibitem{Ladeveze02a}
P.~Ladev{\`e}ze, A.~Nouy, and O.~Loiseau.
\newblock A multiscale computational approach for contact problems.
\newblock {\em Computer Methods in Applied Mechanics and Engineering},
  191(43):4869--4891, 2002.

\bibitem{Ladeveze03a}
Pierre Ladev{\`e}ze and Anthony Nouy.
\newblock On a multiscale computational strategy with time and space
  homogenization for structural mechanics.
\newblock {\em Computer Methods in Applied Mechanics and Engineering},
  192(28-30):3061--3087, 2003.

\bibitem{Mandel93}
Jan Mandel.
\newblock Balancing domain decomposition.
\newblock {\em Communications in Numerical Methods in Engineering}, 9:233--241,
  1993.

\bibitem{Nezamabadi10}
Saeid Nezamabadi, Hamid Zahrouni, Julien Yvonnet, and Michel Potier-Ferry.
\newblock A multiscale finite element approach for buckling analysis of
  elastoplastic long fiber composites.
\newblock {\em International Journal for Multiscale Computational Engineering},
  8(3):287--301, 2010.

\bibitem{Nilsson90}
Karl-Fredrik Nilsson and A~E Giannakopoulos.
\newblock Finite element simulation of delamination growth.
\newblock In M~H Aliabadi, C~A Brebbia, and D~J Cartwright, editors, {\em Proc.
  of the first International Conference on Computer-Aided Assesment and Control
  of Localized Damage}, pages 299--313, Portsmouth, UK, 1990. Springer-Verlag.

\bibitem{Pradeilles-Duval04}
Rachel-Marie Pradeilles~Duval.
\newblock Quasi-static evolution of delaminated structures: analysis of
  stability and bifurcation.
\newblock {\em International Journal of Solids and Structures}, 41(1):103--130,
  2004.

\bibitem{Qiu01}
Y.~Qiu, M.~A. Crisfield, and G.~Alfano.
\newblock An interface element formulation for the simulation of delamination
  with buckling.
\newblock {\em Engineering Fracture Mechanics}, 68(16):1755--1776, 2001.

\bibitem{Storakers88}
Bertil Stor{\aa}kers and B{\"o}rje Andersson.
\newblock Nonlinear plate theory applied to delamination in composites.
\newblock {\em Journal of the Mechanics and Physics of Solids}, 36(6):689--718,
  1988.

\bibitem{Timo61}
S.~P. Timoshenko and J.~M. Gere.
\newblock {\em Theory of Elastic Stability}.
\newblock McGraw-Hill Book Co., New York, 2nd edition, 1961.

\bibitem{Whitcomb89}
J.~D. Whitcomb.
\newblock 3-dimensional analysis of a postbuckled embedded delamination.
\newblock {\em Journal of Composites Materials}, 23(9):862--889, 1989.

\end{thebibliography}

%\begin{thebibliography}{9}
%\end{thebibliography}

\end{document}